\documentclass[10pt]{article}

\usepackage{amsfonts} 
\usepackage{amsmath,hhline,epsfig,latexsym}
\usepackage{amssymb}
\usepackage[latin1]{inputenc}

\usepackage{hyperref}

\usepackage{color}

\newtheorem{theorem}{\sc Theorem}[section]
\newtheorem{lemma}{\sc Lemma}[section]
\newtheorem{proposition}{\sc Proposition}[section]

\newtheorem{remark}{Remark}

\newcommand{\dis}{\displaystyle}

\newcommand{\Fin}{\hfill$\Box$}
\newcommand{\jnt}{\dis\int}
\newcommand{\jjntQT}{\jnt\!\!\!\!\jnt_{Q_{T}}}

\newcommand{\jjntK}{\jnt\!\!\!\!\jnt_{K}}

\newcommand{\N}{\mbox{$I \kern -4pt N$}}
\newcommand{\Q}{\mbox{$Q \kern -8pt I$}}
\newcommand{\R}{\mbox{$I \kern -4pt R$}}
\newcommand{\C}{\mbox{$C \kern -8pt I$}}

\providecommand{\tabularnewline}{\\}

\def\R{{\bf R}}
\title{\textbf{Numerical controllability of the wave equation through primal methods and Carleman estimates}}

\author{
\textsc{Nicolae C\^{\i}ndea}\thanks{Laboratoire de Math\'ematiques, Universit\'e Blaise Pascal (Clermont-Ferand 2), UMR CNRS 6620, Campus de C\'ezeaux, 63177, Aubi\`ere, France. E-mails: {\tt nicolae.cindea@math.univ-bpclermont.fr, arnaud.munch@math.univ-bpclermont.fr.}},
\quad
\textsc{Enrique Fern\'andez-Cara}\thanks{Dpto.~EDAN, Universidad de Sevilla, Aptdo.~1160, 41012, Sevilla, Spain. E-mail: {\tt cara@us.es.}}
\ \ and \ 
\textsc{Arnaud M\"unch}${}^*$}

\date{}

\begin{document}
\maketitle

\begin{abstract}
   This paper deals with the numerical computation of boundary null controls for the 1D wave equation with a potential. The goal is to compute approximations of controls that drive the solution from a prescribed initial state to zero at a large enough controllability time. We do not apply in this work the usual duality arguments but explore instead a direct approach in the framework of global Carleman estimates. More precisely, we consider the control that minimizes over the class of admissible null controls a functional involving weighted integrals of the state and the control. The optimality conditions show that both the optimal control and the associated state are expressed in terms of a new variable, the solution of a fourth-order elliptic problem defined in the space-time domain. We first prove that, for some specific weights determined by the global Carleman inequalities for the wave equation, this problem is well-posed. Then, in the framework of the finite element method, we introduce a family of finite-dimensional approximate control problems and we prove a strong convergence result. Numerical experiments confirm the analysis. We complete our study with several comments.
\end{abstract}

\vskip 0.25cm

\noindent
\textbf{Mathematics Subject Classification (2010)-} 35L10, 65M12, 93B40.



\section{Introduction. The null controllability problem}
\label{sec:intro}

   We are concerned in this work with the null controllability for the 1D wave equation with a potential. The state equation is the following:
\begin{equation}
 \label{eq:wave}
 \left\{
   \begin{array}{ll}
   y_{tt} - (a(x) y_x)_x + b(x, t) y= 0, & \qquad (x, t) \in (0, 1) \times (0, T) \\
   y(0, t) = 0, \qquad y(1, t) = v(t), & \qquad t \in (0, T) \\
   y(x, 0) = y_0(x), \qquad y_t(x, 0) = y_1(x), & \qquad x \in (0, 1).
   \end{array} 
 \right.
\end{equation}
   Here, $T > 0$ and we assume that {\color{black} $a\in C^3([0,1])$ } with $a(x)\geq a_0>0$ in~$[0,1]$, $b \in L^\infty((0, 1) \times (0, T))$, $y_0 \in L^2(0, 1)$ and~$y_1 \in H^{-1}(0, 1)$; $v = v(t)$ is the \textit{control} (a function in~$L^2(0,T)$) and $y=y(x,t)$ is the associated state.

   In the sequel, for any $\tau > 0$ we denote by $Q_\tau$ and $\Sigma_{\tau}$ the sets $(0, 1) \times (0, \tau)$ and $\{0,1\}\times (0,\tau)$, respectively. We will also use the following notation:
   \begin{equation}
\label{eq:L}
L\, y:=y_{tt}-(a(x)y_x)_x + b(x,t)y.
   \end{equation}
 
   For any $(y_0,y_1)\in \boldsymbol{Y}:=L^2(0,1)\times H^{-1}(0,1)$ and any $v\in L^2(0,T)$, it is well known that there exists exactly one solution $y$ to \eqref{eq:wave}, with the following regularity:
\begin{equation}
	\label{eq:wavesol}
	y \in C^0([0, T]; L^2(0, 1)) \cap C^1([0, T]; H^{-1}(0, 1))
\end{equation}
(see for instance~\cite{JLL88}).

   On the other hand, for any $T>0$, the null controllability problem for (\ref{eq:wave}) at time $T$ is the following: for each $(y_0,y_1) \in \boldsymbol{Y}$, find $v\in L^2(0,T)$ such that the corresponding solution to (\ref{eq:wave}) satisfies
   \begin{equation}
\label{eq:nullT}
y(\cdot\,,T) = 0,\ \ y_t(\cdot\,,T) = 0 \quad \text{in }\ (0,1).
   \end{equation}
   
   In view of the linearity and reversibility of the wave equation, \eqref{eq:wave} is null-controllable at $T$ if and only if it is {\it exactly controllable} in~$\boldsymbol{Y}$ at time $T$, i.e.~if and only if for any~$(y_0,y_1) \in \boldsymbol{Y}$ and any~$(z_0,z_1) \in \boldsymbol{Y}$ there exist controls $v\in L^2(0,T)$ such that the associated $y$ satisfies
   $$
y(\cdot\,,T) = z_0,\ \ y_t(\cdot\,,T) = z_1 \quad \text{in }\ (0,1).
   $$

   It is well known that \eqref{eq:wave} is null-controllable at any {\it large} time $T > T^{\star}$ for some $T^{\star}$ that depends on $a$ 
(for instance, see~\cite{BLR, JLL88} for $a \equiv 1$ and $b\equiv 0$ leading to $T^{\star}=2$ and see~\cite{Yao} for a general situation). As a consequence of the {\it Hilbert Uniqueness Method} of J.-L.~Lions~\cite{JLL88}, it is also known that the null controllability of~\eqref{eq:wave} is equivalent to an observability inequality for the associated adjoint problem.

   The goal of this paper is to design and analyze a numerical method allowing to solve the previous null controllability problem.

   So far, the approximation of the minimal $L^2$-norm control --- the so-called HUM control --- has focused most of the attention. The earlier contribution is due to Glowinski and Lions in~\cite{GL96} (see also \cite{Glo08} for an update) and relies on duality arguments. Duality allows to replace the original constrained minimization problem by an unconstrained and \textit{a priori} easier minimization (dual) problem. However, as observed in~\cite{GL96} and later in~\cite{zuazua05}, depending on the approximation method that is used, this approach can lead to some numerical difficulties.
   
   Let us be more precise. It is easily seen that the HUM control is given by $v(t)= a(1)\phi_x(1,t)$, where $\phi$ solves the backwards wave system
   \begin{equation}
\label{heat-phi}
\left\{
\begin{array}{ll}
\dis L \phi = 0 & \quad \textrm{in } Q_T \\
\dis \phi=0 & \quad \textrm{on } \Sigma_T  \\
\dis (\phi(\cdot\,,T),\phi_t(\cdot\,,T))=(\phi_{0}, \phi_{1}) & \quad \textrm{in } (0,1)
\end{array}  
\right.
   \end{equation}
and $(\phi_{0},\phi_{1})$ minimizes the strictly convex and coercive functional 
   \begin{equation}
\mathcal{I}(\phi_0,\phi_1)=\frac{1}{2}\Vert a(1) \phi_x(1,\cdot)\Vert^2_{L^2(0,T)}+\int_0^1 y_0(x)\,\phi_t(x,0)\,dx  -  \langle y_1, \phi(\cdot\,,0)\rangle_{H^{-1},H_0^1}  \label{I_L2norm}
   \end{equation}
over $\boldsymbol{H}=H^1_0(0,1)\times L^2(0,1)$.  Here $\langle \cdot\,,\cdot \rangle_{H^{-1},H_0^1}$ denotes the duality product for~$H^{-1}(0,1)$ and $H_0^1(0,1)$.

   The coercivity of $\mathcal{I}$ over $\boldsymbol{H}$ is a consequence of the {\it observability inequality}
   \begin{equation}
\Vert \phi_0\Vert^2_{H^1_0(0,1)}+ \Vert \phi_1\Vert^2_{L^2(0,1)} \leq C \Vert \phi_x(1,\cdot)\Vert^2_{L^2(0,T)}
\quad \forall (\phi_0,\phi_1)\in \boldsymbol{H},  \label{io}
   \end{equation}
that holds for some constant $C=C(T)$. This inequality has been derived in \cite{JLL88} using the {\it multipliers method.}

   At the numerical level, for standard approximation schemes (based on finite difference or finite element methods), the discrete version of (\ref{io}) may not hold uniformly with respect to the discretization parameter, say $h$. In other words, the constant $C=C(h)$ may blow up as $h$ goes to~zero. Consequently, in such cases the functional $\mathcal{I}_h$ (the discrete version of $\mathcal{I}$) fails to be coercive uniformly with respect to $h$ and the sequence $\{v_h\}_{h>0}$ may not converge to $v$ as $h\to 0$, but diverge exponentially. These pathologies, by now well-known and understood, are due to the spurious discrete high frequencies generated by the finite dimensional approximation; we refer to~\cite{zuazua05} for a review on that topic; see~\cite{munch05} for detailed examples of the behavior observed with finite difference methods.
   
   Several remedies based on more elaborated approximations have been proposed and analyzed in the last decade. Let us mention the use of mixed finite elements~\cite{Castro08}, additional viscosity terms which have the effect to restore the uniform property~\cite{asch, munch05} and also filtering technics~\cite{EZ11}. Also, notice that some error estimates have been obtained recently, see~\cite{Cindea11,EZ11}.   

   In this paper, following the recent work \cite{EFC-AM-11} devoted to the heat equation, we consider a different approach. Specifically, we consider the following extremal problem: 
   \begin{equation}
\label{P-FI}
\left\{
\begin{array}{l}
\dis \hbox{Minimize }\ J(y,v) = {1 \over 2} \jjntQT \rho^2 |y|^2 \,dx\,dt + {1 \over 2} \int_0^T \rho_0^2 |v|^2 \,dt \\
\noalign{\smallskip}
\hbox{Subject to }\ (y,v) \in \mathcal{C}(y_0,y_1;T)
\end{array}
\right.
   \end{equation}
where $\mathcal{C}(y_0,y_1;T)$ denotes the linear manifold
   $$
\mathcal{C}(y_0,y_1;T) = \{\, (y,v) : v \in L^2(0, T),\ \hbox{$y$ solves (\ref{eq:wave})\, and\, satisfies (\ref{eq:nullT})} \,\}.
   $$
   
   Here, we assume that the weights $\rho$ and $\rho_0$ are strictly positive, continuous and uniformly bounded from below by a positive constant in~$Q_T$ and~$(0,T)$, respectively. 
 
   As in the previous $L^2$-norm situation (where we simply have $\rho \equiv 0$ and $\rho_0 \equiv 1$), we can apply duality arguments in order to find a solution to (\ref{P-FI}), by  introducing the unconstrained dual problem
   \begin{equation}
\label{P-FIstar}
\left\{
\begin{array}{l}
\dis \hbox{Minimize }\ J^{\star}(\mu,\phi_0,\phi_1) = {1 \over 2} \jjntQT \rho^{-2} |\mu|^2 \,dx\,dt 
+ {1 \over 2} \jnt_0^T \rho_0^{-2} \vert a(1)\phi_x(1,t)\vert^2 \,dt \\
\noalign{\smallskip}
\dis \phantom{\hbox{Minimize }\ J^{\star}(\mu,\phi_0,\phi_1) } + \jnt_0^1 y_0(x)\,\phi_t(x,0)\,dx  -  \langle y_1,\phi(\cdot\,,0)\rangle_{H^{-1},H_0^1} \\
\noalign{\smallskip}
\hbox{Subject to }\ (\mu,\phi_0,\phi_1) \in L^2(Q_T)\times \boldsymbol{H} ,
\end{array}
\right.
   \end{equation}
   where $\phi$ solves the nonhomogeneous backwards problem 
   \begin{equation}
   \left\{
\begin{array}{ll}
\dis L \phi = \mu & \quad \textrm{in } Q_T \\
\dis \phi=0 & \quad \textrm{on } \Sigma_T  \\
\dis (\phi(\cdot\,,T),\phi_t(\cdot\,,T))=(\phi_{0}, \phi_{1}) & \quad \textrm{in } (0,1).
\end{array}  
   \right.
   \nonumber
   \end{equation}
   
   Here, $J^{\star}$ is the conjugate function of $J$ in the sense of Fenchel and Rockafellar~\cite{temam, Rock} and, if~$\rho\in L^{\infty}(Q_T)$ and~$\rho_0\in L^{\infty}(0,T)$ (that is, $\rho^{-2}$ and $\rho_0^{-2}$ are positively bounded from below), $J^{\star}$ is coercive in $L^2(Q_T)\times \boldsymbol{H}$ thanks to~(\ref{io}). Therefore, if $(\hat{\mu},\hat{\phi_0},\hat{\phi_1})$ denotes the minimizer of $J^{\star}$, the corresponding optimal pair for $J$ is given by
   $$
v=-a(1)\rho_0^{-2} \phi_x(1,\cdot) \ \text{ in } \ (0,T) \ \text{ and } \ y=-\rho^{-2}\mu \ \text{ in } \ Q_T.
   $$
   
   At the discrete level, at least for standard approximation schemes, we may suspect that the coercivity of $J^{\star}$ may not hold uniformly with respect to the discretization parameters, leading to the pathologies and the lack of convergence we have just mentioned.
 
   On the other hand, the fact that the state variable $y$ appears explicitly in the cost $J$ makes it possible to avoid dual methods. We can use instead suitable primal methods to get an optimal pair $(y,v)\in\mathcal{C}(y_0,y_1;T)$.
   The formulation, analysis and practical implementation of these primal methods is the main goal of this paper.
   
   More precisely, the optimality conditions for the functional $J$ allow to express explicitly the optimal pair $(y,v)$ in terms of a new variable, the solution of a fourth-order \textit{elliptic} problem in the space-time domain $Q_T$ that is well-posed under some conditions on $T$, the coefficient $a$ and the weights $\rho$ and $\rho_0$.
   Sufficient conditions are deduced from an appropriate global Carleman estimate, an updated version of the inequalities established in~\cite{Bau01}.
   From a numerical viewpoint, this elliptic formulation is appropriate for a standard finite element analysis. By introducing adequate finite dimensional spaces, we are thus able to deduce satisfactory convergence results for the control, something that does not seem easy to get in the framework of a dual approach.
  
  A similar primal approach, based on ideas by Fursikov and Imanuvilov~\cite{FurIma}, has been used in~\cite{EFC-AM-11} for the numerical null controllability of the heat equation. 

\

   This paper is organized as follows.
   
   In Section~\ref{sec:var}, adapting the arguments and results in~\cite{EFC-AM-11},  we show that the solution to~\eqref{P-FI} can be expressed in terms of the unique solution $p$ to the variational problem (\ref{eq:varp}) in the Hilbert space $P$, defined as the completion of $P_0$ with respect to the inner product (\ref{IP}); see~Proposition~\ref{prop:yp}. The well-posedness is deduced from the application of Riesz's Theorem: a suitable global Carleman inequality ensures the continuity of the linear form in~(\ref{eq:varp}) for $T$ large enough when $\rho$ and $\rho_0$ are given by~(\ref{eq:rho}); see~Theorem~\ref{th:Carleman}.
   
   In Section~\ref{sec:numerical_analysis}, we analyze the variational problem (\ref{eq:varp}) from the viewpoint of the finite element theory. Thus, we replace $P$ by a conformal finite element space $P_h$ of $C^1(\overline{Q_T})$ functions defined by~(\ref{19a}) and we show that the unique solution $\hat{p}_h \in P_h$ to the finite dimensional problem (\ref{pb_phh}) converges (strongly) for the $P$-norm to $p$ as $h$ goes to zero.
   
   Section~\ref{sec:numerics} contains some numerical experiments that illustrate and confirm the convergence of the sequence $\{\hat{p}_h\}$.
   
   Finally, we present some additional comments in~Section~\ref{sec:end_remarks} and we provide some details of the proof of~Theorem~\ref{th:Carleman} in the Appendix. 


\section{A variational approach to the null controllability problem}\label{sec:var}

   With the notation introduced in Section~\ref{sec:intro}, the following result holds:
   
\begin{proposition}\label{prop:solMin}
   Let $T>0$ be large enough. Let us assume that $\rho$ and~$\rho_0$ are positive and satisfy $\rho \in C^0(Q_T)$, $\rho_0 \in C^0(0,T)$ and~$\rho,\rho_0 \geq \underline{\rho} > 0$. Then, for any $(y_0,y_1)\in \boldsymbol{Y}$, there exists exactly one solution to the extremal problem \eqref{P-FI}.  
\end{proposition}

   The proof is simple. Indeed, for $T\geq T^{\star}$, null controllability holds and $\mathcal{C}(y_0,y_1;T)$ is non-empty. Furthermore, it is a closed convex set of $L^2(Q_T)\times L^2(0,T)$. On the other hand, $(y,v) \mapsto J(y,v)$ is strictly convex, proper and lower-semicontinuous in~$L^2(Q_T)\times L^2(0,T)$ and
   $$
J(y,v)\rightarrow +\infty \quad\textrm{as}\quad \Vert (y,v)\Vert_{L^2(Q_T)\times L^2(\Sigma_T)}\rightarrow +\infty.
   $$
Hence, the extremal problem (\ref{P-FI}) certainly possesses a unique solution.

\
  
   In this paper, it will be convenient to assume that the coefficient $a$ belongs to the family
   \begin{equation}\label{eq:setA}
\begin{array}{l}
\dis \mathcal{A}(x_0,a_0) = \{\, {\color{black} a\in C^3([0,1]) } : a(x)\geq a_0 \!>\! 0, \\
\dis \phantom{\mathcal{A}(x_0,a_0) } -\min_{[0,1]} \,\Bigl( a(x) + (x-x_0)a_x(x) \Bigr) < \min_{[0,1]} \,\Bigl( a(x) + {1\over2}(x-x_0)a_x(x) \Bigr) \,\}
\end{array}
   \end{equation}
where $x_0 < 0$ and~$a_0$ is a positive constant.
      
   
   It is easy to check that the constant function $a(x) \equiv a_0$ belongs to $\mathcal{A}(x_0,a_0)$. Similarly, any non-decreasing smooth function bounded from below by~$a_0$ belongs to $\mathcal{A}(x_0,a_0)$. Roughly speaking, $a \in \mathcal{A}(x_0,a_0)$ means that $a$ is sufficiently smooth, strictly positive and not too decreasing in~$[0,1]$.

   Under the assumption \eqref{eq:setA}, there exists ``good" weight functions $\rho$ and $\rho_0$ which provide a very suitable solution to the original null controllability problem. They can be deduced from global Carleman inequalities.

   The argument is the following. First, let us introduce a constant $\beta$, with
   \begin{equation}\label{def_beta}
-\min_{[0,1]} \,\Bigl( a(x) + (x-x_0)a_x(x) \Bigr) < \beta < \min_{[0,1]} \,\Bigl( a(x) + {1\over2}(x-x_0)a_x(x) \Bigr)
   \end{equation}
and let us consider the function
   \begin{equation}
\label{eq:phi}
\phi(x, t):= |x-x_0|^2 - \beta t^2 + M_0,
   \end{equation}
where $M_0$ is such that
   \begin{equation}\label{eq:phiM0}
\phi(x, t) \geq 1 \quad \forall (x, t) \in (0, 1) \times (-T, T),
   \end{equation}
i.e.~$M_0 \geq 1 - |x_0|^2 + \beta T^2$. Then, for any $\lambda > 0$ we set
   \begin{equation}\label{eq:varphi}
\varphi(x, t):= e^{\lambda \phi(x, t)} .
   \end{equation}

   The Carleman estimates for the wave equation are given in the following result:

\begin{theorem}\label{th:Carleman}
   Let us assume that $x_0 < 0$, $a_0 > 0$ and $a \in \mathcal{A}(x_0,a_0)$. Let $\beta$ and $\varphi$ be given respectively by~\eqref{def_beta} and~\eqref{eq:varphi}. Moreover, let us assume that
   \begin{equation}\label{eq:Tlarge}
T > {1 \over \beta} \, \max_{[0, 1]} \,a(x)^{1/2} (x - x_0).
   \end{equation}
Then there exist positive constants $s_0$ and~$M$, only depending on~$x_0$, $a_0$, {\color{black} $\Vert a\Vert_{C^3([0,1])}$, } $\Vert b\Vert_{L^{\infty}(Q_T)}$ and~$T$, such that, for all $s > s_0$, one has
   \begin{equation}\label{eq:CarlT}
   \begin{array}{c}
\dis s \int_{-T}^T \int_0^1 e^{2s\varphi} \left( |w_t|^2 + |w_x|^2 \right) \,dx\,dt +  s^3 \int_{-T}^T  \int_0^1 e^{2s\varphi} |w|^2 \,dx\,dt \\
\noalign{\smallskip}
\dis \leq M \int_{-T}^T \int_0^1 e^{2s\varphi} |L w |^2 \,dx\,dt  + M s \int_{-T}^{T} e^{2s\varphi} |w_x(1,t)|^2 \, dt
   \end{array}
   \end{equation}
for any $w \in L^2(-T, T; H_0^1(0,1))$ satisfying $L w \in L^2((0, 1) \times (-T, T))$ and $w_x(1, \cdot) \in L^2(-T, T)$.
\end{theorem}

   There exists an important literature related to (global) Carleman estimates for the wave equation. Almost all references deal with the particular case $a\equiv 1$; we refer to~\cite{Bau01, Bau11, Ima02, Tataru, Zhang00}. The case where $a$ is non-constant is less studied; we refer to~\cite{FuYongZhang2007}.

   The proof of Theorem~\ref{th:Carleman} follows closely the ideas used in the proofs of~Theorems~2.1 and~2.5 in~\cite{Bau11} to obtain a global Carleman estimate for the wave equation when $a \equiv 1$. The parts of the proof which become different for non-constant $a$ are detailed in the Appendix of this paper.

%


\

   {\color{black}
   In the sequel, it is assumed that $x_0 < 0$ and $a_0 > 0$ are given,  $a \in \mathcal{A}(x_0,a_0)$ and
   \begin{equation}\label{eq:Tlargebis}
T > {2 \over \beta} \, \max_{[0, 1]} \, a(x)^{1/2} (x - x_0), \text{ with } \beta \text{ satisfying } \eqref{def_beta}.
   \end{equation}
   }

   Let us consider the linear space
   $$
P_0 = \{\, q \in C^{\infty}(\overline{Q_T}) : q= 0 \textrm{ on } \Sigma_T \,\}.
   $$ 
 The bilinear form   
   \begin{equation}
\left( p, q \right)_P:= \jjntQT \rho^{-2} Lp \, Lq \,dx\,dt + \int_0^T \rho_0^{-2} \, a(1)^2 \, p_x(1, t) \, q_x(1, t) \,dt 
\label{IP}
   \end{equation}
is a scalar product in~$P_0$. Indeed, in view of~\eqref{eq:Tlargebis}, the unique continuation property for the wave equation holds. Accordingly, if $q \in P_0$, $L q=0$ in~$Q_T$ and~$q_x=0$ on~$\{1\}\times (0,T)$, then $q \equiv 0$. This shows that $(\cdot\,,\cdot)_P$ is certainly a scalar product in~$P_0$.

   Let $P$ be the completion of $P_0$ with respect to this scalar product. Then $P$ is a Hilbert space  for $(\cdot\,,\cdot)_P$ and we can deduce from Theorem~\ref{th:Carleman} the following result, that indicates which are the appropriate weights $\rho$ and $\rho_0$ for our controllability problem:

\begin{lemma}\label{lem:cont}
   Let us assume that $s > s_0$, let us set
   \begin{equation}\label{eq:rho}
\rho(x, t):=  e^{-s\varphi(x, 2t - T)} , \quad \rho_0(t):= \rho(1, t)
   \end{equation}
and let us consider the corresponding Hilbert space $P$. Then there exists a constant $C_0 > 0$, only depending on~$x_0$, $a_0$, $\Vert a\Vert_{C^3([0,1])}$, $\Vert b\Vert_{L^{\infty}(Q_T)}$, $\lambda$, $s$ and~$T$, such that
   \begin{equation}\label{eq:coerc}
\|p(\cdot\,, 0)\|_{H_0^1(0, 1)}^2 + \|p_t(\cdot\,, 0)\|_{L^2(0, 1)}^2 \leq C_0 \left( p, p \right)_P \quad \forall p\in P.
   \end{equation}
\end{lemma}

\noindent
\textsc{Proof:} For every $p \in P$, we denote by $\overline{p} \in L^2((0, 1) \times (-T,T))$ the function defined by
   \[ 
\overline{p}(x,t) = p\left(x, \frac{t+ T}{2}\right).
   \]
It is easy to see that $\overline{p}\in L^2(-T,T;H_0^1(\Omega))$, $L\overline{p} \in L^2((0,1)\times (-T,T))$ and~$\overline{p}_x(1,\cdot)\in L^2(-T,T)$, so that we can apply Theorem~\ref{th:Carleman} to~$\overline{p}$. Accordingly, we have
   \begin{equation}\label{eq:Carlemanp}
   \begin{array}{c}
\dis s \int_{-T}^{T} \int_0^1 e^{2s\varphi} \left( | \overline{p}_t|^2 + |\overline{p}_x|^2 \right) \,dx\,dt  
+ s^3 \int_{-T}^T \int_0^1 e^{2s\varphi}  |\overline{p}|^2 \,dx\,dt \\
\noalign{\smallskip}
\dis \leq C \int_{-T}^T \int_0^1 e^{2s\varphi} |L \overline{p}|^2 \,dx \,dt 
+ Cs \int_{-T}^T e^{2s\varphi(1,t)} |\overline{p}_x(1, t)|^2 \,dt
   \end{array}
   \end{equation}
where $C$ depends on~$x_0$, $a_0$, {\color{black} $\Vert a\Vert_{C^3([0,1])}$, } $\Vert b\Vert_{L^{\infty}(Q_T)}$ and~$T$.

   Replacing $\overline{p}$ by its definition in~\eqref{eq:Carlemanp} and changing the variable $t$ by $t' = 2t-T$ we obtain the following for any $T$ satisfying (\ref{eq:Tlargebis}):
\begin{align*}
  s \jjntQT\rho^{-2} ( | p_t|^2 + |p_x|^2) \,dx \,dt  + s^3 \jjntQT \rho^{-2} |p|^2 \,dx \,dt \\
  \leq C \jjntQT \rho^{-2} |L p|^2 \,dx \,dt + Cs \int_0^T \rho_0^{-2} |p_x(1, t)|^2 \,dt,
\end{align*}
where $C$ is replaced by a slightly different constant. Finally, from Corollary~2.8 in~\cite{Bau11}, we obtain the estimate (\ref{eq:coerc}).
\Fin

\
 
\begin{remark}{\rm
   The estimate (\ref{eq:coerc}) must be viewed as an {\it observability inequality.} As expected, it holds if and only if $T$ is large enough. Notice that, when $a(x) \equiv 1$, the assumption \eqref{eq:Tlargebis} reads
   $$
T > 2(1-x_0) \,.
   $$
This confirms that, in this case, whenever $T > 2$, (\ref{eq:coerc}) holds (it suffices to choose $x_0$ appropriately and apply Lemma~\ref{lem:cont}; see~\cite{JLL88}).
\Fin }
\end{remark}

   The previous results lead to a very useful characterization of the optimal pair $(y,v)$ for $J$:
   
\begin{proposition}\label{prop:yp}
   Let us assume that $s > s_0$, let us set $\rho$ and $\rho_0$ as in~\eqref{eq:rho} and let us consider the corresponding Hilbert space $P$. Let $(y, v) \in \mathcal{C}(y_0, y_1, T)$ be the solution to~\eqref{P-FI}. Then there exists $p \in P$ such that
   \begin{equation}\label{eq:yvp}
y = -\rho^{-2} Lp, \qquad v = - \bigl.(a(x) \rho_0^{-2} p_x)\bigr|_{x=1}.  
   \end{equation}
Moreover, $p$ is the unique solution to the following variational equality:
   \begin{equation}\label{eq:varp}
\left\{
\begin{array}{l}
\dis \jjntQT \rho^{-2} L p\, L q \,dx \,dt  + \int_0^T \rho_0^{-2} a^2(1) p_x(1, t)\, q_x(1, t)\,dt \\
\noalign{\smallskip}
\dis = \int_0^1 y_0(x) \, q_t(x,0) \,dx - \langle y^1, q(\cdot,0)\rangle_{H^{-1}, H_0^1} \quad \forall q \in P; \quad p \in P.
\end{array}
\right.
   \end{equation}
Here and in the sequel, we use the following duality pairing:
   $$
\langle y^1, q(\cdot,0)\rangle_{H^{-1}, H_0^1} = \int_0^1 {\partial \over \partial x} ((-\Delta)^{-1} y_1)(x) \, q_x(x,0)\,dx ,
   $$
where $-\Delta$ is the Dirichlet Laplacian in~$(0,1)$. 
\end{proposition}

\noindent
\textsc{Proof:} From the definition of the scalar product in $P$, we see that $p$ solves \eqref{eq:varp} if and only if
   \[
\left( p, q\right)_P =  \int_0^1 y_0(x) \, q_t(x,0) \,dx - \langle y^1, q(\cdot,0)\rangle_{H^{-1}, H_0^1} \quad \forall q \in P; \quad p \in P.
   \]
 
   In view of Lemma~\ref{lem:cont} and {\it Riesz's Representation Theorem,} problem (\ref{eq:varp}) possesses exactly one solution in $P$.

   Let us now introduce $y$ and $v$ according to (\ref{eq:yvp}) and let us check that $(y,v)$ solves (\ref{P-FI}). First, notice that $y\in L^2(Q_T)$ and $v\in L^2(0,T)$. Then, by replacing $y$ and $v$ in~\eqref{eq:varp}, we obtain the following:
   \begin{equation}\label{eq:varpy}
\jjntQT y\, L q \,dx\,dt + \int_0^T \!\! a(1) v(t) q_x(1, t) \,dt = \int_0^1 \!\! y_0(x) \, q_t(x,0) \,dx - \langle y^1, q(\cdot,0)\rangle_{H^{-1}, H_0^1} \ \ \forall q \in P.
   \end{equation}
Hence, $(y, v)$ is the solution of the controlled wave system \eqref{eq:wave} in the transposition sense. Since $y \in L^2(Q_T)$ and $v \in L^2(0, T)$ the couple $(y, v)$ belongs to $\mathcal{C}(y_0, y_1, T)$. 

   It remains to check that $(y, v)$ minimizes the cost function $J$ in~(\ref{P-FI}). But this is easy. Indeed, for any $(z, w) \in \mathcal{C}(y_0, y_1, T)$ such that $J(z, w) < +\infty$, one has:
   \begin{align*}
 J(z, w) & \geq J(y, v) + \jjntQT \rho^2 y\,(z -y) \,dx \,dt + \int_0^T \rho_0^2 v(w - v) \,dt \\
 & = J(y, v) - \jjntQT Lp\, (z - y) \,dx \,dt + \int_0^T \rho_0^2 v(w - v) \,dt = J(y, v).
   \end{align*}
   The last equality follows from the fact that
   $$
   \begin{aligned}
   \jjntQT Lp\, (z - y) \,dx \,dt = & \jjntQT p\, L(z - y) \,dx \,dt \\
   & + \int_0^1 [p_t\, (z-y)]_0^T \,dx - [<(z-y)_t,p >_{H^{-1},H_0^1}]_0^T \\
   & - \int_0^T [a(x)p_x\,(z-y)]_0^1\,dt  + \int_0^T [a(x) p\, (z-y)_x]_0^1\, dt,
   \end{aligned}
   $$
   the boundary condition for $p$ (see Remark~\ref{order2-4} below), the fact that both $(y,v)$ and $(z,w)$ belong to~$\mathcal{C}(y_0,y_1;T)$ and~(\ref{eq:yvp}).
\Fin

\begin{remark}\label{order2-4}{\rm
   From (\ref{eq:yvp}) and (\ref{eq:varp}), we see that the function $p$ furnished by Proposition~\ref{prop:yp} solves, at least in the distributional sense, the following differential problem, that is of the fourth-order in time and space:
   \begin{equation}
\label{second-fourth}
\left\{
\begin{array}{ll}
\dis L(\rho^{-2} Lp) = 0,                   & (x,t) \in Q_T \\
\dis p(0,t) = 0, \ \  (\rho^{-2} Lp)(0,t)= 0,                                         & t \in (0,T) \\
\dis p(1,t) = 0, \ \  (\rho^{-2} Lp + a\rho_0^{-2}p_x)(1,t)= 0 ,                                        & t \in (0,T) \\
\dis (\rho^{-2} Lp)(x,0) = y_{0}(x), \ \  (\rho^{-2} Lp)(x,T) = 0,   & x \in (0,1) \\
\dis (\rho^{-2} Lp)_t(x,0) = y_1(x), \ \  (\rho^{-2} Lp)_t(x,T) = 0,   & x \in (0,1).
\end{array}
\right.
   \end{equation}
   Notice that the ``boundary'' conditions at $t = 0$ and $t = T$ are of the Neumann kind.
  \Fin}
\end{remark}

\begin{remark}{\rm
   The weights $\rho^{-1}$ and $\rho_0^{-1}$ behave exponentially with respect to $s$. For instance, we have
   $$
\rho(x,t)^{-1}=\exp\,\left( s\, e^{\lambda (\vert x-x_0\vert^2 -\beta (2t-T)^2+M_0)} \right) . 
   $$
For large values of the parameter $s$ (greater than $s_0>0$, see the statement of Theorem \ref{th:Carleman}), the weights $\rho^{-2}$ and $\rho_0^{-2}$ may lead in practice to numerical overflow. One may overcome this situation by introducing a suitable change of variable.

   More precisely, let us introduce the variable $z = \rho p$ and the Hilbert space $M= \rho P$, so that the formulation (\ref{eq:varp}) becomes: 
   \begin{equation}
 \label{eq:varp_m}
 \left\{
  \begin{array}{l}
 \displaystyle \jjntQT \rho^{-2} L (\rho z)\, L(\rho \overline{z}) \,dx\,dt  + \int_0^T \rho_0^{-2} a^2(1) (\rho z)_x(1, t)\, (\rho \overline{z})_x(1, t) \,dt \\
 \displaystyle = \int_0^1 y_0(x)\, (\rho \overline{z})_t(x,0) \,dx - \langle  y_1,(\rho \overline{z})(\cdot,0)\rangle_{H^{-1}, H_0^1}
 \quad \forall \overline{z} \in M; \quad z \in M.
 \end{array}
  \right.
   \end{equation}
   The well-posedness of this formulation is a consequence of the well-posedness of (\ref{eq:varp}). Then, after some computations, the following is found: 
   $$
\begin{aligned}
\rho L(\rho^{-1}z) & = \rho^{-1}\biggl((\rho z)_t - (a(\rho z)_x )_x + b\rho z \biggr) \\
& = (\rho^{-1}\rho_t)z + z_t - a_x ( (\rho^{-1}\rho_x)z+z_x) - a( 2\rho^{-1}\rho_x z_x + \rho^{-1} \rho_{xx} z + z_{xx}) + b\,z
\end{aligned}
   $$
with 
   $$
\rho^{-1}\rho_x = -s \varphi_x(x,2t-T),\quad \rho^{-1}\rho_t = -2s \varphi_t(x,2t-T),\quad\rho^{-1}\rho_{xx} = -s\varphi_{xx} + (s\varphi_x)^2.
   $$
Similarly, 
   $$
(\rho_0^{-1} (\rho z)_x)(1,t)= z_x(1,t).
   $$
Consequently, in the bilinear part of  (\ref{eq:varp_m}), there is no exponential (but only polynomial) function of~$s$. In the right hand side (the linear part), the change of variable introduces negative exponentials in $s$. A similar trick  has been used in~\cite{EFC-AM-11} in the context of the heat equation, where we find weights that blow up exponentially as $t\to T^{-}$. 
\Fin }
\end{remark}

\begin{remark}\label{re:about_s}{\rm
   The exponential form of the weights $\rho$ and $\rho_0$ is purely technical and is related to Carleman estimates. Actually, since for any $s$ and $\lambda$ these weights are uniformly bounded and uniformly positive in~$\overline{Q_T}$, the space $P$ is independent of $\rho$ and $\rho_0$ and one could apply the primal approach to the cost $J$ (defined in (\ref{P-FI})) for any bounded and positive weights. In particular, one could simply take $\rho \equiv 1$ and $\rho_0 \equiv 1$; the estimates (\ref{eq:coerc}) would then read as follows: 
   \begin{equation}
\|p(\cdot\,, 0)\|_{H_0^1(0, 1)}^2 + \|p_t(\cdot\,, 0)\|_{L^2(0, 1)}^2 \leq C_0 \biggl( \Vert L\,p \Vert^2_{L^2(Q_T)} + \Vert a(1)\,p_x(1,\cdot)\Vert^2_{L^2(0,T)} \biggr) \quad \forall p\in P
   \end{equation}
for some constant $C_0>0$. 
   This inequality can also be obtained directly by the {\it multipliers method;} we refer to~\cite{Yao} and references therein.
\Fin }
\end{remark}

\begin{remark}\label{re:newrho0}{\rm
  As remarked in \cite{Bau11} (see Remark 2.7), the estimate~\eqref{eq:coerc} can be proven for a weight $\rho_0$ which blows up at $t = 0$ and $t = T$. For this purpose, we consider a function $\theta_\delta \in C^2([0, T])$ with $\theta_\delta(0) = \theta_\delta(1) = 0$ and $\theta_\delta(x) = 1$ for every $x \in (\delta, T-\delta)$. Then, introducing again $\overline{p}(x,t):= \theta_\delta(t)p(x, (t+T)/2)$, it is not difficult to see that the proofs of Lemma~\ref{lem:cont} and~Theorem~\ref{th:Carleman} can be adapted to obtain \eqref{eq:coerc} with
   \[
\rho(x, t) =  e^{-s\varphi(x, 2t - T)} , \qquad \rho_0(t) = \theta_\delta(t)^{-1}\rho(1, t).
   \]
Thanks to the properties of~$\theta_\delta$, the control $v$ defined by
   $$
v = -\theta_{\delta}^2 \bigl. \rho_0^{-2} a(1)p_x \bigr|_{x=1}
   $$
vanishes at~$t=0$ and also at~$t=T$, a property which is very natural and useful in the boundary controllability context. In the sequel, we will use this modified weight $\rho_0$, imposing in addition, for numerical purposes, the following behavior near $t=0$ and $t=T$: 
   \begin{equation}
\lim_{t\rightarrow 0^{+}}  \frac{\theta_{\delta}(t)}{\sqrt{t}}=\mathcal{O}(1), \quad  \lim_{t\rightarrow T^{-}}  \frac{\theta_{\delta}(t)}{\sqrt{T-t}}=\mathcal{O}(1).  \label{hyp_theta}
   \end{equation}
\Fin }
\end{remark}


\section{Numerical analysis of the variational approach}
\label{sec:numerical_analysis}

We now highlight that the variational formulation (\ref{eq:varp}) allows to obtain a sequence of  approximations $\{v_h\}$ that converge strongly towards the null control $v$ furnished by the solution to~\eqref{P-FI}.


\subsection{A conformal finite dimensional approximation}

Let us introduce the bilinear form $m(\cdot,\cdot)$ over $P\times P$
   $$
m(p,q) := (p,q)_P = \jjntQT \rho^{-2}Lp \, Lq \,dx\,dt + \int_0^T  a(1)^2 \rho_0^{-2} p_x(1,t) \, q_x(1,t) \,dt
   $$
and the linear form $\ell$, with  
   $$
\langle \ell,q \rangle:=   \jnt_0^1 y_0(x) \, q_t(x,0) \,dx - \langle y^1, q(\cdot,0)\rangle_{H^{-1}, H_0^1}.
   $$ 
Then (\ref{eq:varp}) reads as follows: 
   \begin{equation}
m(p,q)=\langle \ell,q \rangle, \quad \forall q\in P; \quad p\in P.  \label{pb_p}
   \end{equation}

   Let us assume that a finite dimensional space $P_h\subset P$ is given for each $h \in \mathbb{R}_+^2$.
   Then we can introduce the following {\it approximated} problems: 
   \begin{equation}
m(p_h,q_h)=\langle \ell,q_h \rangle, \quad \forall q_h\in P_h; \quad p_h\in P_h.   \label{pb_ph}
   \end{equation}
   Obviously, each (\ref{pb_ph}) is well-posed. Furthermore, we have the following classical result: 

\begin{lemma}\label{lemma1}
   Let $p\in P$ be the unique solution to \eqref{pb_p} and let $p_h \in P_h$ be the unique solution to \eqref{pb_ph}. Then we have:
   \begin{equation}
\Vert p-p_h\Vert_P \leq \inf_{q_h\in P_h} \Vert p - q_h \Vert_P. \nonumber
   \end{equation}
\end{lemma}

\noindent
\textsc{Proof:}
   We write that
   $$
\Vert p_h-p\Vert_P^2= m(p_h-p,p_h-p) = m(p_h-p,p_h-q_h) + m(p_h-p,q_h-p).
   $$
   The first term vanishes for all $q_h\in P_h$.
   The second one is bounded by $\Vert p_h-p\Vert_P \Vert q_h-p\Vert_P$.
   So, we get 
   $$
\Vert p-p_h \Vert_P \leq \Vert p - q_h\Vert_P \quad \forall q_h\in P_h
   $$
and the result follows. \Fin

\

   As usual, this result can be used to prove that $p_h$ converges towards $p$ when the spaces $P_h$ are chosen appropriately. 
   More precisely, let us assume that an {\it interpolation operator} $\Pi_h: P_0\to P_h$ is given for any $h\in \mathbb{R}_+^2$ and let us suppose that
   \begin{equation}
\Vert p - \Pi_h p \Vert_P\to 0 \ \textrm{ as } \ h\to (0,0) \quad  \forall p\in P_0. \label{convergence_H1}
   \end{equation} 
   We then have the following convergence result:
   
\begin{proposition}\label{convergence_zhz}
   Let $p\in P$ be the solution to~\eqref{pb_p} and let $p_h\in P_h$ be the solution to~\eqref{pb_ph} for each $h \in \mathbb{R}_+^2$.
   Then
   \begin{equation}\label{30p}
\Vert p-p_h\Vert_P \rightarrow 0 \ \textrm{ as }\ h \to (0,0).
   \end{equation}
\end{proposition}

\noindent
\textsc{Proof:}
   Let us choose $\epsilon>0$.
   Since $P_0$ is dense in~$P$, there exists $p_\epsilon\in P_0$ such that $\Vert p - p_\epsilon \Vert_P\leq \epsilon$.
   Therefore, we find from Lemma~\ref{lemma1} that
   \begin{equation}
\begin{aligned}
\Vert p-p_h\Vert_P  & \leq \Vert p- \Pi_h p_\epsilon\Vert_P  \\
&  \leq \Vert p-p_\epsilon \Vert_P + \Vert p_\epsilon-\Pi_h p_\epsilon\Vert_P \\
& \leq \epsilon + \Vert p_\epsilon-\Pi_h p_\epsilon\Vert_P.
\end{aligned}
\nonumber
\end{equation} 
But we know from (\ref{convergence_H1}) that $\Vert p_\epsilon-\Pi_h p_\epsilon\Vert_P$ goes to zero as $h \in \mathbb{R}_+^2$, $h\to (0,0)$. Consequently, we also have \eqref{30p}.
\Fin


\subsection{The finite dimensional spaces $P_h$}

   The spaces $P_h$ must be chosen such that $\rho^{-1}Lp_h$ belongs to $L^2(Q_T)$ for any $p_h\in P_h$.
   This means that $p_h$ must possess second-order derivatives in~$L^2_{\rm loc}(Q_T)$.
   Therefore, a conformal approximation based on a standard quadrangulation of $Q_T$ requires spaces of functions continuously differentiable with respect to both variables $x$ and~$t$.

   For large integers $N_{x}$ and $N_{t}$, we set $\Delta x = 1/N_{x}$, $\Delta t = T/N_{t}$ and $h = (\Delta x,\Delta t)$.
   We introduce the associated quadrangulations ${\mathcal Q}_{h}$, with $\overline{Q_T} = \bigcup_{K \in {\mathcal Q}_{h}} K$ and we assume that  $\{\mathcal Q_{h}\}_{h>0}$ is a regular family.
   Then, we introduce the space $P_h$ as follows:
   \begin{equation}
   \label{19a}
   P_{h} = \{\, z_{h} \in C^1(\overline{Q_T}) : z_h|_{K} \in  \mathbb{P}(K) \ \ \forall K \in  {\mathcal Q}_{h}, \ z_h = 0 \ \text{ on } \Sigma_T \,\} .
   \end{equation}
   
   Here, $\mathbb{P}(K)$ denotes the following space of polynomial functions in $x$ and $t$:
  \begin{equation}
\mathbb{P}(K)=(\mathbb{P}_{3,x} \otimes \mathbb{P}_{3,t})(K) \label{def_PK}
  \end{equation}
where $\mathbb{P}_{r,\xi}$ is by definition the space of polynomial functions of order $r$ in the variable $\xi$.

   Obviously, $P_h$ is a finite dimensional subspace of $P$.


   {\color{black}
   Let us introduce the notation
   $$
K_{kl}:= [x_k,x_{k+1}]\times [t_l,t_{l+1}],
   $$
where
   $$
x_k:=(k-1)\Delta x, \ t_l:=(l-1)\Delta t, \ \text{ for } \ k=1,\dots,N_x+1, \ l=1,\dots,N_t+1.
   $$
  
   For any $k$, we denote by $(L_{ik})_{0\leq i\leq 3}$ the Hermite functions associated to~$[x_k,x_{k+1}]$. They are given by 
   $$
\left\{
\begin{aligned}
& L_{0k}(x):=(1+2c)(1-c)^2, \quad L_{1k}(x):=c^2(3-2c) \\
& L_{2k}(x):=\Delta x\, c(1-c)^2, \quad L_{3k}(x):= \Delta x\, c^2(c-1) \\
&c:=(x-x_k)/\Delta x.
\end{aligned}
\right.
   $$
   Recall that, for any $f\in C^1([x_k,x_{k+1}])$, the function
   $$
(\Pi_{\Delta x}f)(x):= \sum_{i=0}^1 L_{ik}(x)f(x_{i+k})+ \sum_{i=0}^1 L_{i+2,k}(x) f_x(x_{i+k}) 
   $$
is the unique element in $\mathbb{P}_3([x_k,x_{k+1}])$ that satisfies
   $$
(\Pi_{\Delta x}f)(x_{k+i})=f(x_{k+i}), \quad (\Pi_{\Delta x}f)_x(x_{k+i})=(f_x)(x_{k+i}), \quad i=0,1.
   $$
   
   In a similar way, we denote by $(\mathcal{L}_{jl})_{0\leq j\leq 3}$ the Hermite functions associated to the time interval~$[t_l,t_{l+1}]$. Then, from the definition of $\mathbb{P}(K_{kl})$, we can obtain easily for any $u\in P_0$ the polynomial function in~$\mathbb{P}(K_{kl})$ uniquely determined by the values of $u$, $u_x$, $u_t$ and~$u_{xt}$ at the vertices of~$K_{kl}$:
   
\begin{lemma}
   For each $u \in P_0$, let us define the function $\Pi_h u$ as follows: for any $k$ and~$l$, 
   \begin{equation}
\begin{aligned}
\Pi_h u(x,t):=& \sum_{i,j=0}^1 L_{ik}(x)\mathcal{L}_{jl}(t) u(x_{i+k},t_{j+l})  + \sum_{i,j=0}^1 L_{i+2,k}(x)\mathcal{L}_{jl}(t) u_x(x_{i+k},t_{j+l}) \\
& + \sum_{i,j=0}^1 L_{ik}(x)\mathcal{L}_{j+2,l}(t) u_t(x_{i+k},t_{j+l}) + \sum_{i,j=0}^1 L_{i+2,k}(x)\mathcal{L}_{j+2,l}(t) u_{xt}(x_{i+k},t_{j+l})
\end{aligned}
\nonumber
   \end{equation}
in~$K_{kl}=[x_k,x_k+\Delta x]\times [t_l,t_l+\Delta t]$.

   Then $\Pi_h u$ is the unique function in $P_h$ that satisfies
   \begin{equation}
\begin{aligned}
\Pi_h u(x_{k+i},t_{l+j}) = u(x_{k+i},t_{l+j}),  \quad (\Pi_h u(x_{k+i},t_{l+j}))_x = u_x(x_{k+i},t_{l+j}),  \\
(\Pi_h u(x_{k+i},t_{l+j}))_t = u_t(x_{k+i},t_{l+j}),  \quad (\Pi_h u(x_{k+i},t_{l+j}))_{xt} = u_{xt}(x_{k+i},t_{l+j})
\end{aligned}
\nonumber
   \end{equation}
for all $i,j\in \{0,1\}$. The linear mapping $\Pi_h : P_0 \mapsto P_h$ is by definition the interpolation operator associated to~$P_h$.
\end{lemma}

}

This result allows to get an expression of  $u-\Pi_h u$ on each element $K_{kl}$ that will be used in the next section:
   
\begin{lemma} \label{diff_u_pu}
   For any $u \in P_0$, we have
   \begin{equation}
   \begin{aligned}
u-\Pi_h u = & \sum_{i,j=0}^1 \biggl( m_{ij} u_x(x_{i+k},t_{j+l}) +  n_{ij} u_t(x_{i+k},t_{j+l}) + p_{ij} u_{tx}(x_{i+k},t_{j+l})\biggr) \\
  &+ \sum_{i,j=0}^1 L_{ik}\mathcal{L}_{jl} \mathcal{R}[u;x_{i+k},t_{j+l}]
\end{aligned}
\label{diff_u_pu_eq}
   \end{equation}
in~$K_{kl}$, where the $m_{ij}$, $n_{ij}$ and~$p_{ij}$ are given by
   \[
   \left\{
   \begin{aligned}
& m_{ij}(x,t):= \biggl(L_{ik}(x)(x-x_i)-L_{i+2,k}(x)\biggr)\mathcal{L}_j(t) \\
& n_{ij}(x,t):= L_{ik}(x)\biggl(\mathcal{L}_j(t) (t-t_j)-\mathcal{L}_{j+2}(t)\biggr) \\
& p_{ij}(x,t):= L_{ik}(x)\mathcal{L}_{jl}(t) (x-x_i)(t-t_j) - L_{i+2}(x)\mathcal{L}_{j+2}(t)
\end{aligned}
\right.
   \]
and
   \[
\begin{array}{l}
\dis \mathcal{R}[u;x_{i+k},t_{j+l}](x,t):= \int_{t_{j+l}}^t (t-s)u_{tt}(x_{i+k},s) \,ds
\\ \noalign{\smallskip} \dis
\phantom{\mathcal{R}[u;x_{i+k},t_{j+l}](x,t):=}+ (x-x_{i+k})\int_{t_{j+l}}^t (t-s)u_{xtt}(x_{i+k},s) \,ds 
\\ \noalign{\smallskip} \dis
\phantom{\mathcal{R}[u;x_{i+k},t_{j+l}](x,t):=} + \int_{x_{i+k}}^x (x-s)u_{xx}(s,t) \,ds.
\end{array}
   \]
\end{lemma}

  The proof is very simple. In fact, \eqref{diff_u_pu_eq} is a consequence of the following Taylor expansion for $u$ with integral remainder: 
   \begin{equation}
\begin{aligned}
u(x,t) &= u(x_k,t_l)+ (t-t_l)u_t(x_k,t_l)+ \int_{t_l}^t (t-s) u_{tt}(x_k,s) \,ds \\
         & \ \ + (x-x_k)\biggl( u_x(x_k,t_l)+(t-t_l)u_{xt}(x_k,t_l)+ \int_{t_l}^t (t-s)u_{xtt}(x_k,s) \,ds\biggr) \\
         & \ \ + \int_{x_k}^x (x-s)u_{xx}(s,t) \,ds
\end{aligned}
\nonumber
   \end{equation}
and the identity $\sum_{i,j=0}^1 L_{ik}(x)\mathcal{L}_{jl}(t) \equiv 1$.


\subsection{An estimate of $\Vert p - \Pi_h p\Vert_P$ and some consequences}

   Let us now prove that (\ref{convergence_H1}) holds when the $P_h$ are given by~\eqref{19a}--\eqref{def_PK}.
 
   Thus, let us fix $p \in P_0$ and let us first check that
   \begin{equation}\label{41p}
\jjntQT \rho^{-2} \vert L (p-\Pi_h(p)) \vert^2\,dx\,dt \to 0 \ \text{ as } \ h=(\Delta x, \Delta t) \to (0,0).
   \end{equation}
   
   For each $K_{kl} \in \mathcal{Q}_h$
   (simply denoted by $K$ in the sequel), we write:
   \begin{equation}
   \label{41pbis}
   \begin{aligned}
\jjntK & \rho^{-2} \vert L(p-\Pi_h p)\vert^2 \,dx\,dt
\leq  \Vert \rho^{-2}\Vert_{L^{\infty}(K)} \jjntK \vert L(p-\Pi_h p)\vert^2 \,dx\,dt \\
& \leq 3 \Vert \rho^{-2}\Vert_{L^{\infty}(K)}\biggl( \jjntK \vert (p-\Pi_h p)_{tt}\vert^2 \,dx\,dt + \jjntK \vert (a(x)(p-\Pi_h p)_x)_x\vert^2 \,dx\,dt \\
& \quad + \Vert b\Vert^2_{L^{\infty}(K)} \jjntK \vert p-\Pi_h p\vert^2 \,dx\,dt \biggr).
\end{aligned}
   \end{equation}
   Using Lemma~\ref{diff_u_pu},  we have:
   \begin{equation}
\begin{aligned}
\jjntK &\vert p -\Pi_h p \vert^2   \, dx\,dt \\
& = \jjntK \biggl|\sum_{i,j}  \bigl( m_{ij} p_x(x_i,t_j) + n_{ij}p_t(x_i,t_j)+p_{ij} p_{tx}(x_i,t_j) + L_i \mathcal{L}_j \mathcal{R}[p;x_{i},t_{j}]\bigr) \biggr|^2\,dx\,dt \\
& \leq 16 \Vert p_x\Vert^2_{L^{\infty}(K)} \sum_{i,j} \jjntK |m_{ij}|^2 \,dx\,dt+ 16 \Vert p_t\Vert^2_{L^{\infty}(K)} \sum_{i,j} \jjntK |n_{ij}|^2 \,dx\,dt \\
& \quad + 16 \Vert p_{tx}\Vert^2_{L^{\infty}(K)} \sum_{i,j} \jjntK |p_{ij}|^2 \,dx\,dt+ 16\sum_{i,j}\jjntK |L_i \mathcal{L}_j \mathcal{R}[p;x_i,t_j]|^2 \,dx\,dt,
\end{aligned}
\nonumber
   \end{equation}
where we have omitted the indices $k$ and~$l$.  

   Moreover, 
   \begin{equation}
\begin{aligned}
|\mathcal{R}[p;x_{i+k},t_{j+l}]|^2 \leq &  \vert t-t_j\vert^3\Vert p_{tt}(x_i,\cdot)\Vert^2_{L^2(t_l,t_{l+1})} + |x-x_i|^2 \vert t-t_j\vert^3  \Vert p_{xtt} (x_i,\cdot)\Vert^2_{L^2(t_l,t_{l+1})}\\
& + \vert x-x_i\vert^3 \Vert p_{xx}(\cdot\,,t)\Vert^2_{L^2(x_k,x_{k+1})}. 
\end{aligned}
\nonumber
   \end{equation}
   Consequently, we get:
   \begin{equation}
\begin{aligned}
\sum_{i,j} & \jjntK |L_i \mathcal{L}_j \mathcal{R}[p;x_{i+k},t_{j+l}]|^2 \,dx\,dt \\
& \leq \sup_{x\in (x_k,x_{k+1})} \Vert p_{tt}(x,\cdot)\Vert^2_{L^2(t_l,t_{l+1})}
\sum_{i,j}\jjntK |L_i(x) \mathcal{L}_j(t)|^2 \vert t-t_j\vert^3 \,dx\,dt \\
& \quad + \sup_{x\in (x_k,x_{k+1})} \Vert p_{xtt}(x,\cdot)\Vert^2_{L^2(t_l,t_{l+1})}
\sum_{i,j}\jjntK |L_i(x) \mathcal{L}_j(t)|^2 \vert t-t_j\vert^3 |x-x_i|^2 \,dx\,dt \\
& \quad + \sup_{t\in (t_l,t_{l+1})} \Vert p_{xx}(\cdot\,,t)\Vert^2_{L^2(x_k,x_{k+1})}
\sum_{i,j}\jjntK |L_i(x) \mathcal{L}_j(t)|^2 \vert x-x_i\vert^3 \,dx\,dt. 
\end{aligned}
\nonumber
   \end{equation}
   After some tedious computations, one finds that
   \begin{equation}
   \begin{aligned}
& \sum_{i,j} \jjntK |m_{ij}|^2 \,dx\,dt = \frac{104}{11025} (\Delta x)^3 \Delta t, \quad \sum_{i,j} \jjntK |n_{ij}|^2 \,dx\,dt = \frac{104}{11025} \Delta x(\Delta t)^3,\\
& \sum_{i,j} \jjntK |p_{ij}|^2 \,dx\,dt = \frac{353}{198450} (\Delta x)^3 (\Delta t)^3
\end{aligned}
\nonumber
   \end{equation}
   and
\begin{equation}
\begin{aligned}
& \sum_{i,j} \jjntK |L_i(x) \mathcal{L}_j(t)\vert^2 \vert t-t_j\vert^3\,dx\,dt= \frac{143}{7350} \Delta x (\Delta t)^4, \\
&  \sum_{i,j} \jjntK |L_i(x) \mathcal{L}_j(t)\vert^2 \vert x-x_i\vert^3\,dx\,dt= \frac{143}{7350} (\Delta x)^4 \Delta t, \\
&  \sum_{i,j} \jjntK |L_i(x) \mathcal{L}_j(t)\vert^2 \vert x-x_i\vert^2 \vert t-t_j\vert^3 \,dx\,dt= \frac{209}{132300} (\Delta x)^3 (\Delta t)^4.\\
\end{aligned}
\nonumber
\end{equation}
This leads to the following estimate for any $K=K_{kl}\in \mathcal{Q}_h$:
   \begin{equation}
\begin{aligned}
\jjntK |p -\Pi_h p|^2 \,dx\,dt  \leq &\frac{1664}{11025} (\Delta x)^3\Delta t   \Vert p_x \Vert^2_{L^{\infty}(K)}  \\
& + \frac{1664}{11025} \Delta x (\Delta t)^3 \Vert p_t \Vert^2_{L^{\infty}(K)} \\
& + \frac{2824}{99225}  (\Delta x)^3 (\Delta t)^3 \Vert p_{tx} \Vert^2_{L^{\infty}(K)} \\
& + \frac{1144}{3675} (\Delta x)^4 \Delta t \sup_{x\in (x_k,x_{k+1}) } \Vert p_{tt}(\cdot\,,t)\Vert_{L^2(t_l,t_{l+1})}^2\\
& + \frac{1144}{3675} \Delta x (\Delta t)^4 \sup_{t\in (t_l,t_{l+1}) } \Vert p_{xx}(\cdot\,,t)\Vert_{L^2(x_k,x_{k+1})}^2\\
& + \frac{836}{33075} (\Delta x)^3 (\Delta t)^3 \sup_{x\in (x_k,x_{k+1}) } \Vert p_{xtt}(\cdot\,,t)\Vert_{L^2(t_l,t_{l+1})}^2.
\end{aligned}
\nonumber
\end{equation}
   We deduce that
      \begin{equation}
\begin{aligned}
\jjntQT |p -\Pi_h p|^2 \,dx\,dt  \leq & K_1 T \Vert p_x \Vert^2_{L^{\infty}(Q_T)}\,  (\Delta x)^2  \\
& + K_1 T \Vert p_t \Vert^2_{L^{\infty}(Q_T)}\, (\Delta t)^2 \\
& + K_2  T  \Vert p_{tx} \Vert^2_{L^{\infty}(Q_T)}\,(\Delta x)^2 (\Delta t)^2 \\
& + K_3 \Vert p_{tt}(\cdot\,,t)\Vert_{L^2(0,T; L^{\infty}(0,1))}^2\, (\Delta x)^3 \\
& + K_3 \Vert p_{xx}(\cdot\,,t)\Vert_{L^{\infty}(0,T; L^2(0,1)}^2\,  (\Delta t)^3 \\
& + K_4 \Vert p_{xtt}(\cdot\,,t)\Vert_{L^2(0,T; L^{\infty}(0,1))}^2\, (\Delta x)^2 (\Delta t)^2
\end{aligned}
\nonumber
\end{equation}
for some positive constants $K_i$. Hence, for any $p\in P_0$ one has
   $$
\jjntQT \vert p-\Pi_h \,p\vert^2\,dx\,dt  \to 0 \ \textrm{ as } \ h \to (0,0).
   $$
   
   Proceeding as above, we show that the other terms in (\ref{41pbis}) also converge to $0$. Hence, (\ref{41p}) holds.
   
   On the other hand, a similar argument yields
   $$
\int_0^T  \rho_0^{-2} a(1)^2 \vert (p-\Pi_h p)_x\vert^2 \,dx\,dt \to 0 \ \textrm{ as } \ h \to (0,0)
   $$
and, consequently, we find that (\ref{convergence_H1}) holds.
   
   We can now use Proposition~\ref{convergence_zhz} and deduce convergence results for the approximate control and state variables:
   
\begin{proposition}\label{propconv-vh-a}
   Let $p_h\in P_h$ be the unique solution to~\eqref{pb_ph}, where $P_h$ is given by~\eqref{19a}--\eqref{def_PK}. Let us set
   $$
y_h:= \rho^{-2} L p_h, \quad v_h:= - \rho_0^{-2} \bigl. a(x) p_{h,x} \bigr|_{x=1}.
   $$
   Then one has
   $$
\| y - y_h \|_{L^2(Q_T)} \to 0 \ \text{and } \ \| v - v_h \|_{L^2(0,T)} \to 0,
   $$
where $(y,v)$ is the solution to~\eqref{P-FI}.
\Fin
\end{proposition}


\subsection{A second approximated problem}

   For simplicity, we will assume in this section that $y_1 \in C^0([0,1])$.
   
   In order to take into account the numerical approximation of the weights and the data that we necessarily have to perform in practice, we will also consider a second approximated problem.
   It is the following:
   \begin{equation}
m_h(\hat{p}_h,q_h)= \langle \ell_h,q_h \rangle \quad \forall q_h\in P_h; \quad \hat{p}_h\in P_h, \label{pb_phh}
   \end{equation} 
where the bilinear form $m_h(\cdot\,,\cdot)$ is given by
   $$
m_h(p_h,q_h):=\jjntQT \pi_h(\rho^{-2}) Lp_h \, Lq_h \,dx\,dt + \int_0^T  a(1)^2 \pi_{\Delta x}(\rho_{0}^{-2}) p_h \, q_h \,dt
   $$
and the linear form $\ell_h$ is given by
   $$
\langle \ell_h,q \rangle:=   \jnt_0^1 (\pi_{\Delta x}y_0)(x) \, q_t(x,0) \,dx - \langle \pi_{\Delta x}y_1, q(\cdot,0)\rangle_{H^{-1}, H_0^1}.
   $$
   
   Here, for any function $f\in C^0(\overline{Q_T})$, $\pi_h(f)$ denotes the piecewise linear function which coincides with $f$ at all vertices of $\mathcal{Q}_h$. Similar (self-explanatory) meanings can be assigned to $\pi_{\Delta x}(z)$ and $\pi_{\Delta t}(w)$ when $z \in C^0([0,1])$ and $w \in C^0([0,T])$, respectively.

   Since the weight $\rho^{-2}$ is strictly positive and bounded in~$Q_T$ (actually~$\rho^{-2} \geq 1$), we easily see that the ratio $\pi_h(\rho^{-2})/ \rho^{-2}$ is bounded uniformly with respect to $h$ (for $|h|$ small enough).  The same holds for the vanishing weight $\theta_{\delta}^2\rho(1,\cdot)^{-2}$ under the assumptions (\ref{hyp_theta}).

   As a consequence, it is not difficult to prove that (\ref{pb_phh}) is well-posed. Moreover, we have:

\begin{lemma}
   Let $p_h$ and $\hat{p}_h$ be the solutions to~\eqref{pb_ph} and~\eqref{pb_phh}, respectively. Then,  
   \begin{equation}\label{phathtop}
\begin{aligned}
\Vert \hat{p}_h -p_h\Vert_P \leq & \max\biggl( \Vert \frac{\pi_h(\rho^{-2})}{\rho^{-2}}-1\Vert_{L^{\infty}(Q_T)}, \Vert \frac{\pi_{\Delta t}(\rho_0^{-2})}{\rho_0^{-2}}-1\Vert_{L^{\infty}(0,T)}\biggr) \Vert \hat{p}_h\Vert_P \\
& + C_1 \Vert \pi_{\Delta x}(y_0)-y_0\Vert_{L^2}  +C_2 \Vert \pi_{\Delta x}(y_1)-y_1\Vert_{H^{-1}},
\end{aligned}
   \end{equation}
where $C_1$ and $C_2$ are positive constants independent of $h$.
\end{lemma}

\noindent
{\sc Proof: }
   Since $p_h$ and $\hat{p}_h$ respectively solve~\eqref{pb_ph} and~\eqref{pb_phh}, one has:
   \begin{equation}
\begin{aligned}
\Vert &\hat{p}_h -p_h\Vert_P^2 =  m(\hat{p}_h-p_h,\hat{p}_h-p_h) \\ 
& = m(\hat{p}_h,\hat{p}_h-p_h) - m_h(\hat{p}_h,\hat{p}_h-p_h) + \langle l_h, \hat{p}_h-p_h\rangle - \langle l, \hat{p}_h-p_h\rangle \\
& = \jjntQT  (\rho^{-2}\!-\!\pi_h(\rho^{-2})) L\hat{p}_h  \, L(\hat{p}_h\!-\!p_h) \,dx\,dt + \int_0^T (\rho_0^{-2}\!-\!\pi_{\Delta t}(\rho_0^{-2})) a(1)^2 \hat{p}_{x,h} 
(\hat{p}_{x,h}\!-\!p_{x,h}) \,dt \\
& \quad + \langle l_h, \hat{p}_h-p_h\rangle - \langle l, \hat{p}_h-p_h\rangle \\
& = \jjntQT  \biggl(1-\frac{\pi_{h}(\rho^{-2})}{\rho^{-2}}\biggr) (\rho^{-1}L\hat{p}_h)  (\rho^{-1}L(\hat{p}_h-p_h)) \,dx\,dt \\
& \quad + \int_0^T \biggl(1-\frac{\pi_{\Delta t}\rho^{-2}_0}{\rho_0^{-2}}\biggr) a(1)^2 (\rho_0^{-1} \hat{p}_h) \rho_0^{-1}(\hat{p}_{h,x}-p_{h,x})\,dt \\
& \quad + \int_{\Omega} (\pi_{\Delta x}(y_0)-y_0)(x) \, (\hat{p}_{t,h}-p_{t,h})(x,0)\,dx - \langle \pi_{\Delta x}(y_1)-y_1, (\hat{p}_{h}-p_{h})(x,0)\rangle_{H^{-1},H_0^1}\,.
\end{aligned}\nonumber
   \end{equation}
In view of the definitions of the bilinear forms $m(\cdot\,,\cdot)$ and $m_h(\cdot\,,\cdot)$, we easily find \eqref{phathtop}.
\Fin

\

   Taking into account that (\ref{convergence_H1}) holds and
   $$
\max\biggl( \Vert \frac{\pi_h(\rho^{-2})}{\rho^{-2}}-1\Vert_{L^{\infty}(Q_T)}, 
\Vert \frac{\pi_{\Delta t}(\rho_0^{-2})}{\rho_0^{-2}}-1\Vert_{L^{\infty}(0,T)}\biggr)
\to 0,
   $$
we find that, as $h$ goes to zero, the unique solution to~(\ref{pb_phh}), converges in~$P$ to the unique solution to~(\ref{pb_p}):
   $$
\Vert p-\hat{p}_h\Vert_P \leq  \Vert p-p_h\Vert_P + \Vert p_h-\hat{p}_h\Vert_P  \rightarrow 0.
   $$
   
   An obvious consequence is the following:
   
\begin{proposition}\label{propconv-vh}
   Let $\hat p_h\in P_h$ be the unique solution to~\eqref{pb_phh}, where $P_h$ is given by~\eqref{19a}--\eqref{def_PK}. Let us set
   \begin{equation}\label{haty-hatv}
\hat y_h:= \pi_{h}(\rho^{-2}) L \hat{p}_h, \quad \hat v_h:= - \pi_{\Delta x}(\rho_0^{-2}) \bigl. a(x) \hat{p}_{h,x} \bigr|_{x=1}.
   \end{equation}
   Then one has
   $$
\| y - \hat y_h \|_{L^2(Q_T)} \to 0 \ \text{and } \ \| v - \hat v_h \|_{L^2(0,T)} \to 0,
   $$
where $(y,v)$ is the solution to~\eqref{P-FI}.
\Fin
\end{proposition}


\section{Numerical experiments}
\label{sec:numerics}

   We now present some numerical experiments concerning the solution of (\ref{pb_phh}), which can in fact be viewed as a linear system involving a banded sparse, definite positive, symmetric matrix of order $4N_x N_t$. We will denote by $\mathcal{M}_h$ this matrix. If $\{\hat{p}_h\}$ stands for the corresponding vector solution of size $4N_xN_t$, we may write  $(\hat{p}_h,q_h)_{P_h}= (\mathcal{M}_h \{\hat{p}_h\},\{q_h\})$ for any $q_h\in P_h$. 

   We will use an exact integration method in order to compute the components of $\mathcal{M}_h$ and the (direct) Cholesky method with reordering to solve the linear system.

   After the computation of $\hat{p}_h$, the control $\hat{v}_h$ is given by (\ref{haty-hatv}). Observe that, in view of the definition of the space $P_h$,  the derivative with respect to $x$ of $\hat{p}_h$ is a degree of freedom of $\{\hat{p}_h\}$; hence, the computation of $\hat{v}_h$ does not require any additional calculus.
   
   The corresponding controlled state $\hat{y}_h$ may be obtained by using the pointwise first equality \eqref{haty-hatv} or, equivalently, by solving (\ref{eq:varpy}). However, in order to check the action of the control function $\hat{v}_h$ properly, we have computed $\hat{y}_h$ by solving (\ref{eq:wave}) with a $C^1$ finite element method in space and a standard centered scheme of second order in time.
   
   Thus, let us introduce the finite dimensional spaces 
   $$
Z_h=\{\, z_h \in C^1([0,1]) : \bigl. z_h \bigr|_{[x_i,x_i+\Delta x]} \in \mathbb{P}_{3,x} \ \ \forall i=1,\dots, N_x \,\}
   $$
and $Z_{0h}=\{\, z_h\in Z_h : z_h(0) = z_h(1) = 0 \,\}$. Then, a suitable approximation $\hat{y}_h$ of the controlled state $y$ is defined in the following standard way:

\begin{itemize}

\item At time $t=0$, $\hat{y}_h$ is given by $y_h(\cdot,0)=P_{Z_h}(y_0)$, the projection of $y_0$ on $Z_h$;

\item At time $t_1=\Delta t$, $\hat{y}_h(\cdot,t_1)\in Z_h$ is given by the solution to
   \begin{equation}\label{post_y1}
\left\{
\begin{array}{l}
\dis 2 \int_0^1  \frac{(\hat{y}_h(x,t_1)- \hat{y}_h(x,t_0)-\Delta t \,y_1(x))}{(\Delta t)^2} \phi\, dx  \\
\noalign{\smallskip}
\dis \quad + \int_0^1 \left[ a(x)\hat{y}_{h,x}(x,t_0) \phi_x +b(x,t_0)\hat{y}_{h}(x,t_0)\phi \right]\,dx=0 \\
\noalign{\smallskip}
\dis \forall \phi \in Z_{h0}; \quad \hat{y}_{h}(0,t_1) \in Z_h, \ \hat{y}_{h}(0,t_1)=0, \ \hat{y}_{h}(1,t_1)=\hat{v}_h(t_1).
\end{array}
\right.
   \end{equation}
   
\item At time $t=t_n= n \Delta t$, $n=2,\cdots, N_t$, $\hat{y}_h(\cdot,t_n)$ solves the following linear problem:
   \begin{equation}\label{post_y2}
\left\{
\begin{array}{l}
\dis 2 \int_0^1  \frac{(\hat{y}_h(x,t_n)- 2\hat{y}_h(x,t_{n-1})+\hat{y}_h(\cdot,t_{n-2})}{(\Delta t)^2} \phi\, dx  \\
\noalign{\smallskip}
\dis \quad + \int_0^1 \left[ a(x)\hat{y}_{h,x}(x,t_{n-1}) \phi_x +b(x,t_n)\hat{y}_{h}(x,t_{n-1})\phi \right]\,dx=0 \\
\noalign{\smallskip}
\dis \forall \phi \in Z_{h0}; \quad \hat{y}_{h}(0,t_n) \in Z_h, \ \hat{y}_{h}(0,t_n)=0, \ \hat{y}_{h}(1,t_n)=\hat{v}_h(t_n).
\end{array}
\right.
   \end{equation}
   
\end{itemize}

   This requires a preliminary projection of $\hat{v}_h$ on a grid on $(0,T)$ fine enough in order to fulfill the underlying CFL condition. To this end, we use the following interpolation formula: for any $p_h\in P_h$ and any $\theta\in [0,1]$, we have:
   \begin{equation}
\label{interpol}
\begin{aligned}
p_{h,x}(1,t_j+\theta \Delta t) &= (2\theta+1)(\theta-1)^2\,p_{h,x}(1,t_j) + \Delta t\,\theta(1-\theta)^2\,p_{h,xt}(1,t_j)  \\
& \quad + \theta^2(3-2\theta)\,p_{h,x}(1,t_{j+1}) + \Delta t\,\theta^2(\theta-1)\,p_{h,xt}(1,t_{j+1})
\end{aligned}
   \end{equation}
for all $t\in [t_j,t_{j+1}]$.

\ 

We will consider a constant coefficient $a(x) \equiv a_0=1$ and a constant potential $b(x,t) \equiv 1$ in $Q_T$. We will take $T=2.2$, $x_0=-1/20$, $\beta=0.99$ and $M_0=1-x_0^2+\beta T^2$, so that (\ref{eq:Tlargebis}) holds. Finally, concerning the parameters $\lambda$ and $s$ (which appear in (\ref{eq:Carlemanp})), we will take $\lambda=0.1$ and $s=1$. 

\begin{remark}{\rm
   Let us emphasize that our approach does not require in any way the discretization meshes to be uniform.
\Fin }
\end{remark}


\subsection{Estimating the Carleman constant}\label{carlemanconstant}

   Before prescribing the initial data, let us check that the finite dimensional analog of the observability constant $C_0$ in~(\ref{eq:coerc}) is uniformly bounded with respect to $h$ when (\ref{eq:Tlargebis}) is satisfied. We consider here the case $a\equiv 1$ and $b\equiv 1$.

   In the space $P_h$, the approximate version of~(\ref{eq:coerc}) is
   $$
(A_h \{p_h\},\{p_h\}) \leq C_{0h} (\mathcal{M}_h \{p_h\},\{p_h\}) \quad  \forall \{p_h\} \in P_h,
   $$
where $A_h$ is the square matrix of order $4N_tN_x$ defined by the identities
   $$
(A_h \{p_h\},\{q_h\}):= \jnt_0^1 \left( p_{h,x}(x,0) \, q_{h,x}(x,0) + p_{h,t}(x,0) \, q_{h,t}(x,0) \right) \,dx.
   $$
   Therefore, $C_{0h}$ is the solution of a generalized eigenvalue problem: 
   \begin{equation}\label{43p}
C_{0h} = \max \{\, \lambda : \exists p_h\in P_h, \ \ p_h \not= 0, \ \text{ such that } \ A_h \{p_h\} = \lambda \mathcal{M}_h \{p_h\} \,\}.
   \end{equation}

   We can easily solve \eqref{43p} by the power iteration algorithm. Table~\ref{tab:C0} collects the values of $C_{0h}$ for various $h = (\Delta x,\Delta t)$ for $T=2.2$ and $T=1.5$, with $\Delta t=\Delta x$. As expected, $C_{0h}$ is bounded in the first case only. The same results are obtained for $\Delta t\neq \Delta x$.

\begin{table}[http]
\centering
\begin{tabular}{|c|cccc|}
\hline
$\Delta x,\Delta t$  & $1/10$ & $1/20$ & $1/40$ & $1/80$  \tabularnewline

\hline $T=2.2$  & $6.60\times 10^{-2}$ & $7.61\times 10^{-2}$ & $8.56\times 10^{-2}$ & $9.05\times 10^{-2}$  \tabularnewline

$T=1.5$  & $0.565$ & $2.672$ & $17.02$ & $96.02$  \tabularnewline

\hline
\end{tabular}
\caption{The constant $C_{0h}$ with respect to $h$.}
\label{tab:C0}
\end{table}

{\color{black}
  In agreement with Remark \ref{re:about_s}, we obtain the same behavior of the constant with respect to $T$ for any $s$, in particular for $s=0$ leading to $\rho\equiv 1$ and $\rho_0\equiv 1$.
  }


\subsection{Smooth initial data and constant speed of propagation}
   We now solve (\ref{eq:varp}) with $a\equiv 1$ and smooth initial data.  For simplicity, we also take a constant potential $b\equiv  1$.

   For $(y_0,y_1)=(\sin(\pi x),0)$, Table~\ref{tab:ex1} collects relevant numerical values with respect to $h=(\Delta x,\Delta t)$. We have taken $\Delta t=\Delta x$ for simplicity but, in this finite element framework, any other choice is possible. In particular, we have reported the condition number $\kappa(\mathcal{M}_h)$ of the matrix $\mathcal{M}_h$, defined by
   $$
\kappa(\mathcal{M}_h) = \vert\vert\vert \mathcal{M}_h\vert\vert\vert_2 \, \vert\vert\vert \mathcal{M}^{-1}_h\vert\vert\vert_2
   $$
(the norm $\vert\vert\vert \mathcal{M}_h\vert\vert\vert_2$ stands for the largest singular value of $\mathcal{M}_h$). We observe that this number behaves polynomially with respect to $h$.

\begin{table}[http]
\centering
\begin{tabular}{|c|cccc|}
\hline
$\Delta x,\Delta t$  & $1/10$ & $1/20$ & $1/40$ & $1/80$  \tabularnewline

\hline $\kappa(\mathcal{M}_h)$  & $3.06\times 10^{8}$ & $1.57\times 10^{10}$ & $6.10\times 10^{11}$ & $2.47\times 10^{13}$  \tabularnewline

$\Vert \hat{p}_h\Vert_{P}$ & $1.541\times 10^{-1}$ & $1.548\times 10^{-1}$ & $1.550\times 10^{-1}$ & $1.550\times 10^{-1}$ \tabularnewline

$\Vert \hat{p}_h-p\Vert_{P}$ & $4.46\times 10^{-2}$ & $1.45\times 10^{-2}$ & $4.01\times 10^{-3}$ & $8.38\times 10^{-4}$ \tabularnewline

$\Vert \hat{v}_h\Vert_{L^2(0,T)}$  & $5.421\times 10^{-1}$ & $5.431\times 10^{-1}$ & $5.434\times 10^{-1}$ & $5.434\times 10^{-3}$ \tabularnewline

$\Vert \hat{v}_{h}-v\Vert_{L^2(0,T)}$  & $2.39\times 10^{-2}$ & $8.12 \times 10^{-3}$ & $2.48 \times 10^{-3}$ & $9.57\times 10^{-4}$ 
\tabularnewline

$\Vert \hat{y}_h(\cdot\,,T)\Vert_{L^2(0,1)}$  & $1.80\times 10^{-2}$ & $8.18 \times 10^{-3}$ & $1.64 \times 10^{-3}$ & $5.85 \times 10^{-4}$ \tabularnewline

$\Vert \hat{y}_{t,h}(\cdot\,,T)\Vert_{H^{-1}(0,1)}$  & $3.06\times 10^{-2}$ & $8.25 \times 10^{-3}$ & $3.59 \times 10^{-3}$ & $1.93 \times 10^{-3}$ 
\tabularnewline\hline
\end{tabular}
\caption{$(y_0(x),y_1(x))\equiv(\sin(\pi x),0)$, $a\equiv 1, b\equiv 1$ - $T=2.2$.}
\label{tab:ex1}
\end{table}

   Table~\ref{tab:ex1} clearly exhibits the convergence of the variables $\hat{p}_h$ and $\hat{v}_h$ as $h$ goes to zero. Assuming that $h=(1/160,1/160)$ provides a reference solution, we have also reported in Table~\ref{tab:ex1} the estimates $\Vert p-\hat{p}_h\Vert_{P}$ and $\Vert v-\hat{v}_h\Vert_{L^2(0,T)}$. We observe then that
   $$
\Vert p-\hat{p}_h\Vert_{P}=\mathcal{O}(h^{1.91}), \quad \Vert v-\hat{v}_h\Vert_{L^2(0,T)}=\mathcal{O}(h^{1.56}).
   $$
The corresponding state $\hat{y}_h$ is computed from the main equation $(\ref{eq:wave})$, as explained above, taking $\Delta t=\Delta x/4$. That is, we use (\ref{interpol}) with $\theta=0,1/4,1/2$ and $3/4$ on each interval $[t_j,t_{j+1}]$. We observe the following behavior with respect to $h$: 
   $$
\Vert \hat{y}_h(\cdot,T)\Vert_{L^2(0,1)}=\mathcal{O}(h^{1.71}), \quad \Vert \hat{y}_{t,h}(\cdot,T)\Vert_{H^{-1}(0,T)}=\mathcal{O}(h^{1.31}),
   $$
which shows that the control $\hat{v}_h$ given by the second equality in~\eqref{haty-hatv} is a good approximation of a null control for (\ref{eq:wave}).

\begin{figure}[http]
\begin{center}
\begin{minipage}[t]{6.9cm}
\hspace*{-1.cm}
\includegraphics[scale=0.44]{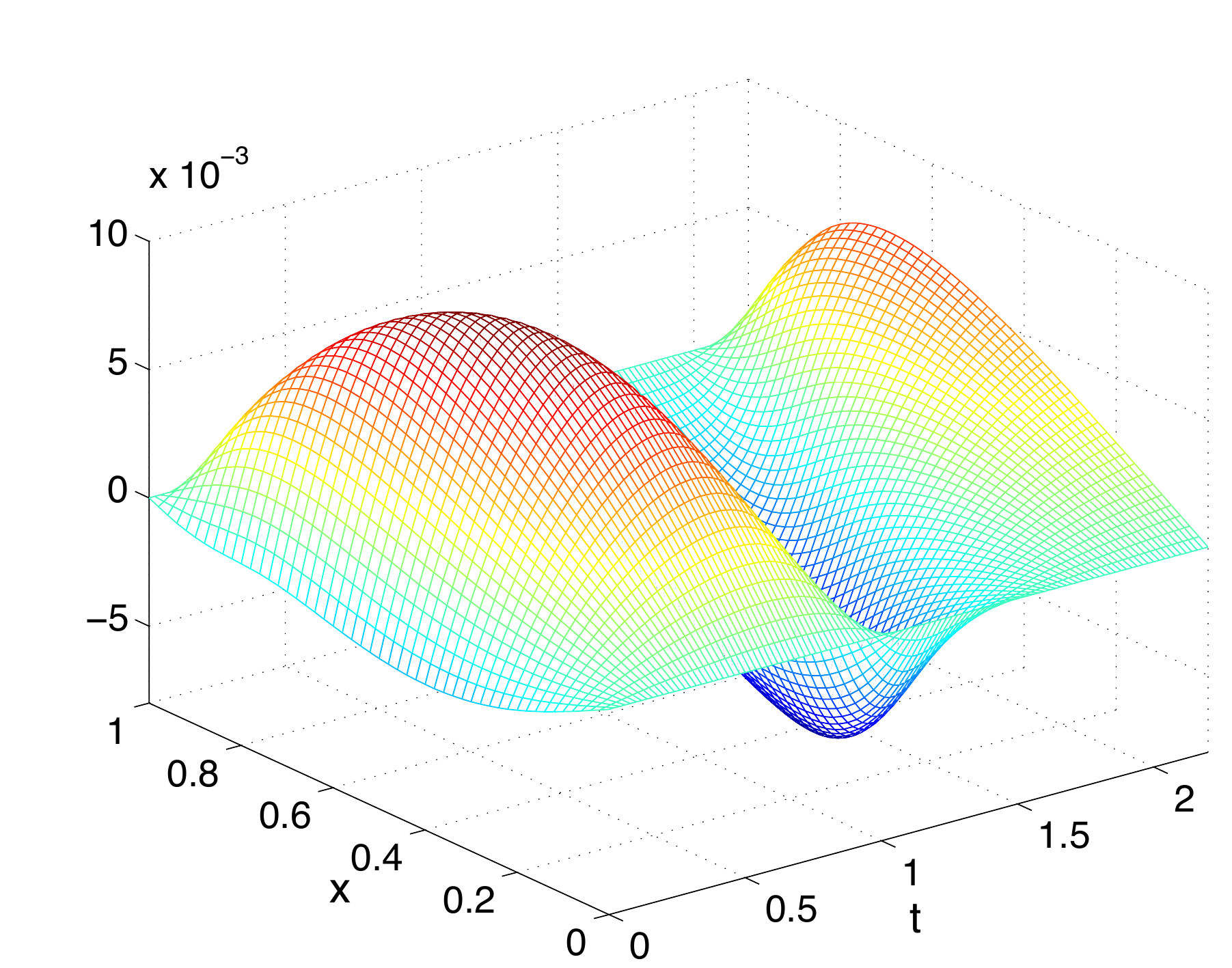}
\end{minipage}
\begin{minipage}[t]{6.9cm}
\hspace*{0.5cm}
\includegraphics[scale=0.38]{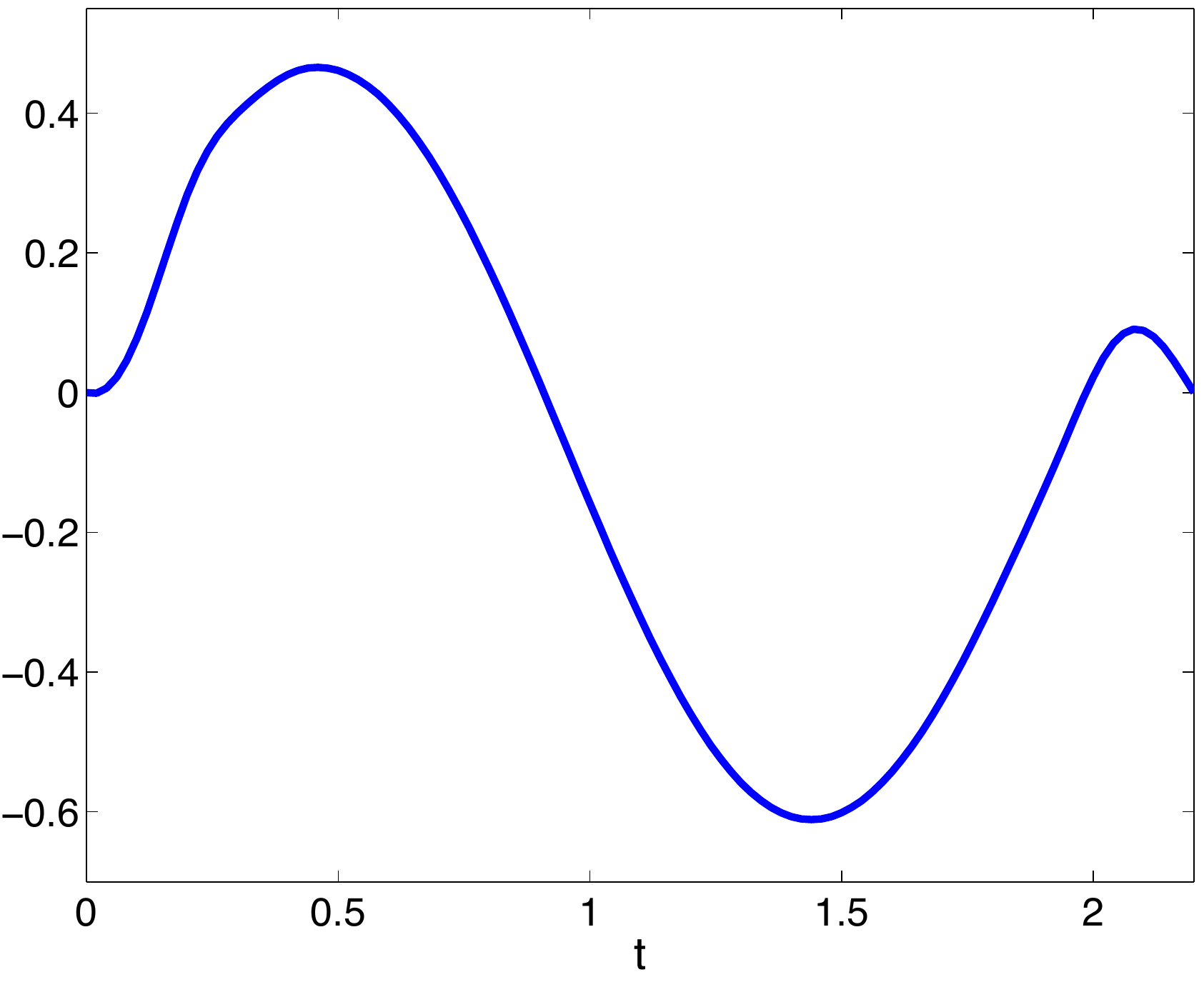}
\end{minipage}
\caption{$(y_0(x),y_1(x))\equiv(\sin(\pi x),0)$ and $a\equiv 1$ - The solution $\hat{p}_h$  over $Q_T$ (\textbf{Left}) and the corresponding variable $\hat{v}_h$ on $(0,T)$ (\textbf{Right}) - $h=(1/80,1/80)$.}\label{fig:p_controle_ex1}
\end{center}
\end{figure}

\begin{figure}[http]
\begin{center}
\includegraphics[scale=0.44]{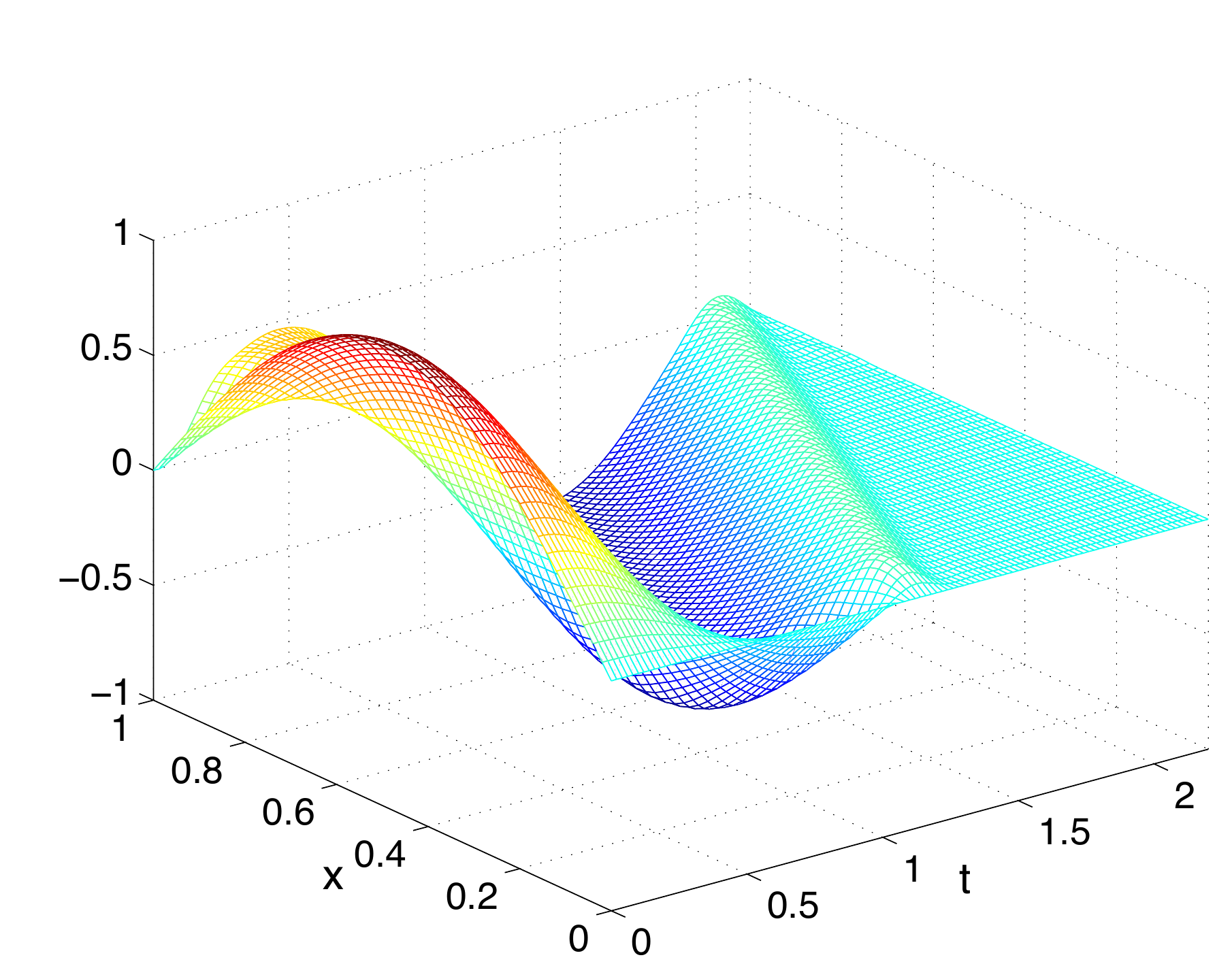}
\caption{$(y_0(x),y_1(x))\equiv(\sin(\pi x),0)$ and $a\equiv 1$ -  The solution $\hat{y}_h$  over $Q_T$ - $h=(1/80,1/80)$.}\label{fig:y_ex1}
\end{center}
\end{figure}

   Figure~\ref{fig:p_controle_ex1}-Left displays the function $\hat{p}_h \in P$ (the unique solution to~(\ref{pb_phh})) for $h=(1/80,1/80)$. Figure~\ref{fig:p_controle_ex1}-Right displays the associated control $\hat{v}_h$. As a consequence of the introduction of the function $\theta_\delta$ in the weight, we see that $\hat{v}_h$ vanishes at times $t=0$ and $t=T$. Finally, Figure~\ref{fig:y_ex1} displays the corresponding controlled state $\hat{y}_h$.

\begin{figure}[http]
\begin{center}
\includegraphics[scale=0.44]{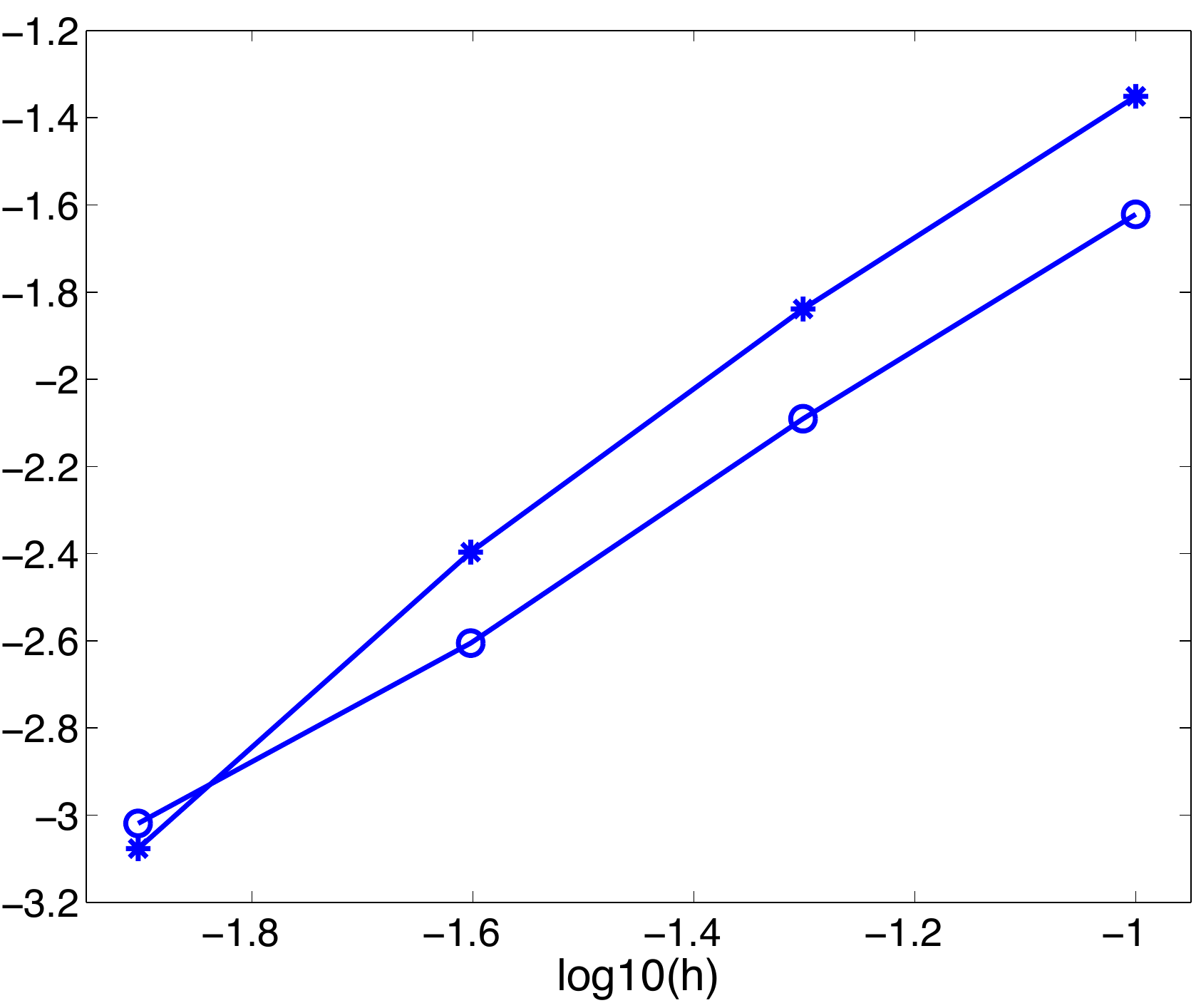}
\caption{$\log_{10}\Vert p-\hat{p}_h\Vert_P$ ($\star$) and $\log_{10}\Vert v-\hat{v}_h\Vert_{L^2(0,T)}$ ($\circ$) vs. $\log_{10}(h)$.}\label{fig:logp_logv_ex1}
\end{center}
\end{figure}

\begin{figure}[http]
\begin{center}
\begin{minipage}[t]{6.9cm}
\hspace*{-1.cm}
\includegraphics[scale=0.44]{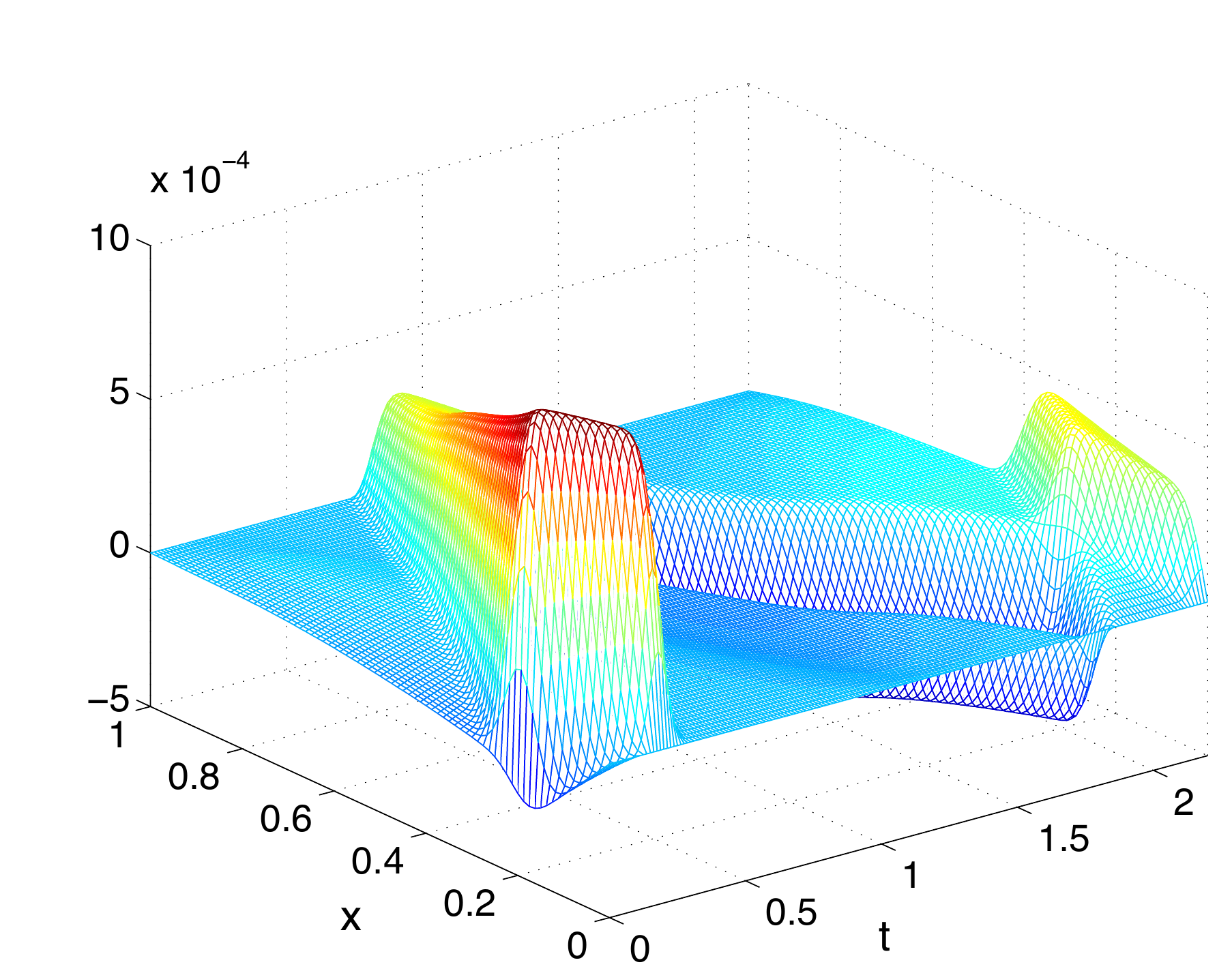}
\end{minipage}
\begin{minipage}[t]{6.9cm}
\hspace*{0.5cm}
\includegraphics[scale=0.38]{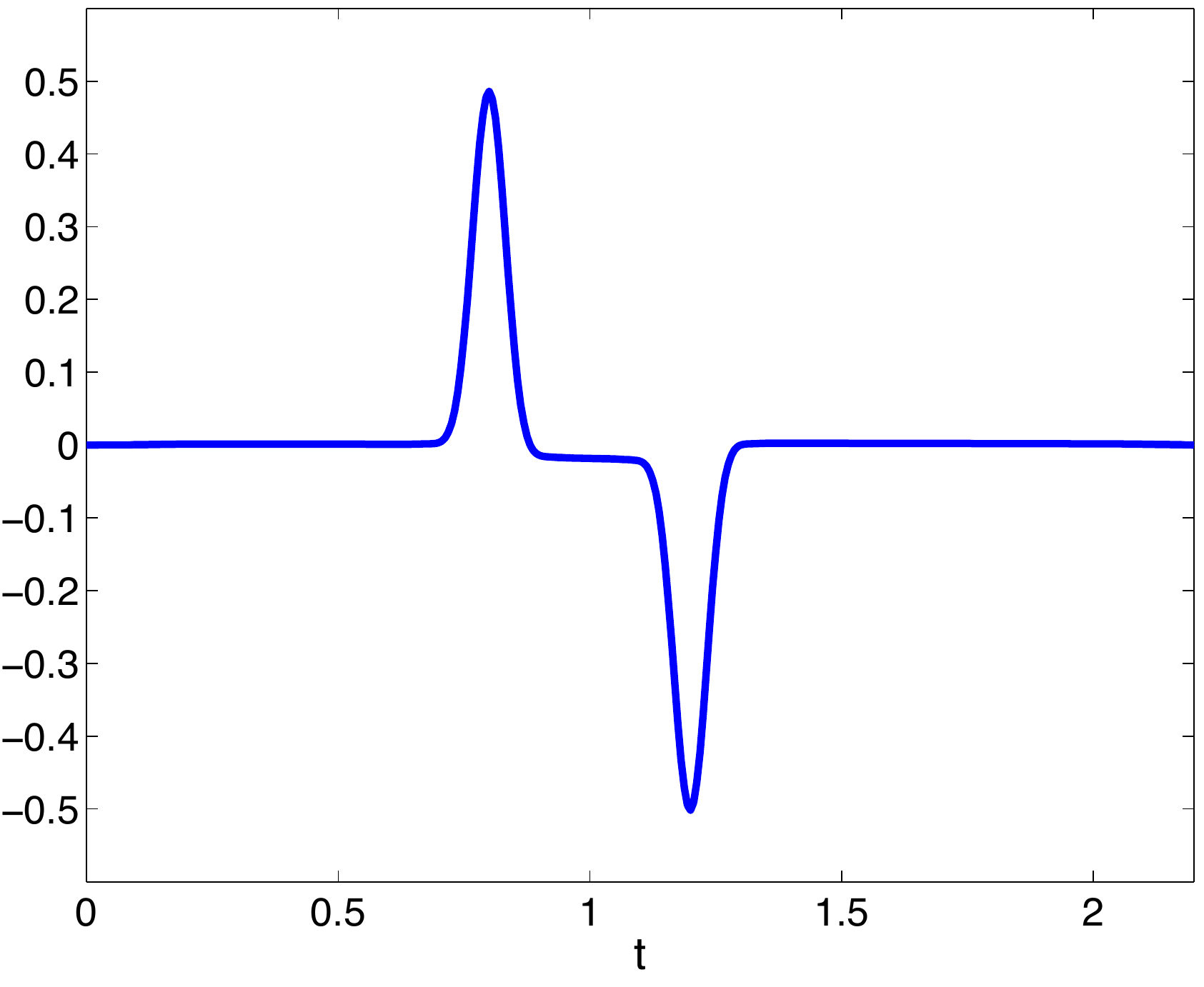}
\end{minipage}
\caption{$(y_0(x),y_1(x))\equiv (e^{-500(x-0.2)^2},0)$ and $a\equiv 1$ - The solution $\hat{p}_h$  over $Q_T$ (\textbf{Left}) and the corresponding variable $\hat{v}_h$ on $(0,T)$ (\textbf{Right}) - $h=(1/80,1/80)$.}\label{fig:p_controle_ex2}
\end{center}
\end{figure}

\begin{figure}[http!]
\begin{center}
\includegraphics[scale=0.5]{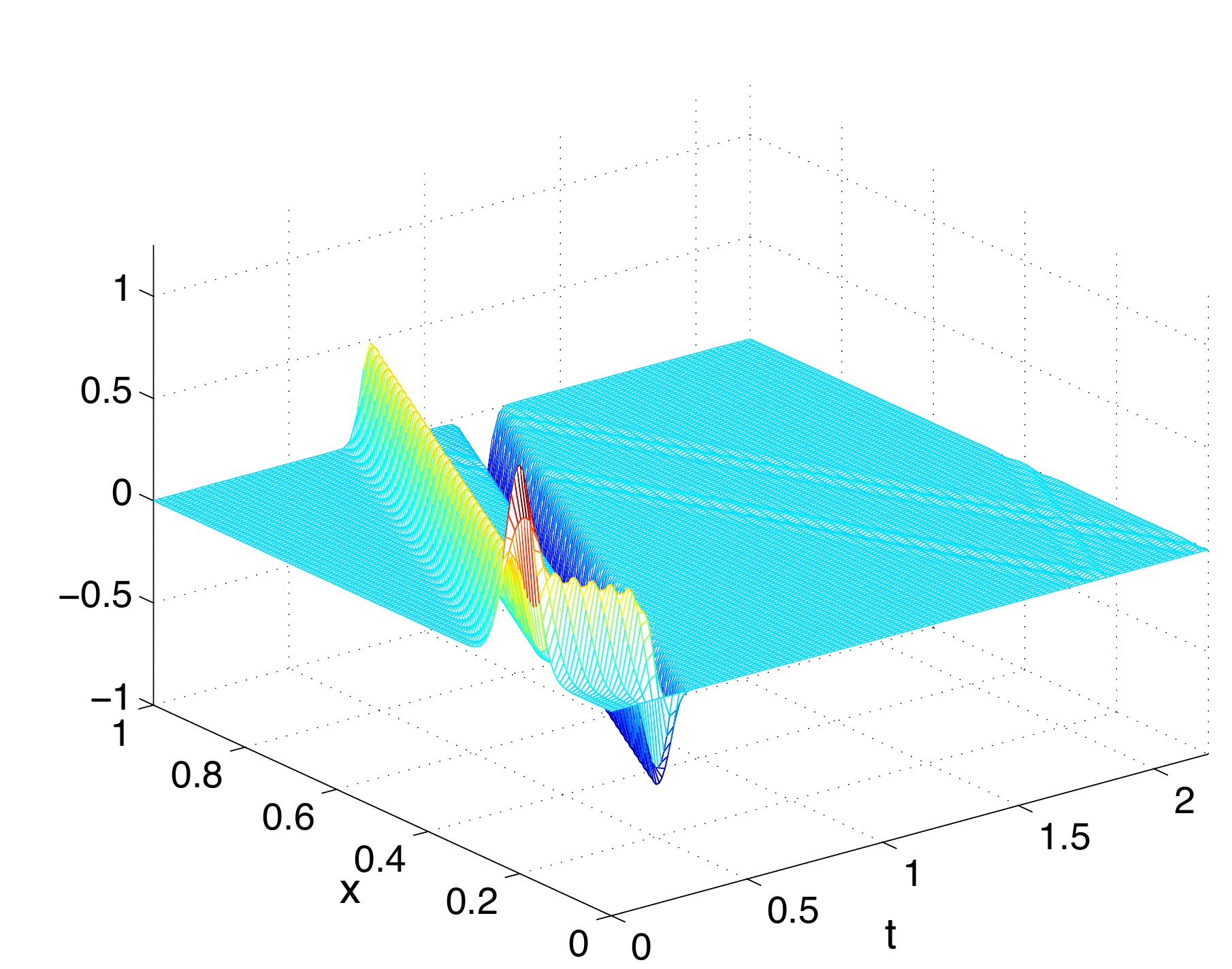}
\caption{$(y_0(x),y_1(x))\equiv (e^{-500(x-0.2)^2},0)$ and $a\equiv 1$ -The solution $\hat{y}_h$  over $Q_T$ - $h=(1/80,1/80)$.}\label{fig:y_ex2}
\end{center}
\end{figure}

   Table~\ref{tab:ex2} and~Figures~\ref{fig:p_controle_ex2} and~\ref{fig:y_ex2} provide the results for $y_0(x) \equiv e^{-500(x-0.2)^2}$ and $y_1(x) \equiv 0$. We still observe the convergence of the variables $\hat{p}_h$, $\hat{v}_h$ and $\hat{y}_h$, with a lower rate. This is due in part to the shape of the initial condition $y_0$. Precisely, we get $\Vert p-\hat{p}_h\Vert_{P} = \mathcal{O}(h^{1.74})$, $\Vert \hat{v}_{h}-v\Vert_{L^2(0,T)}=\mathcal{O}(h^{0.68})$, $\Vert \hat{y}_h(\cdot,T)\Vert_{L^2(0,1)}=\mathcal{O}(h^{1.35})$ and $\Vert \hat{y}_{t,h}(\cdot,T)\Vert_{H^{-1}(0,T)}=\mathcal{O}(h^{1.11})$.
   
\begin{table}[http]
\centering
\begin{tabular}{|c|ccccc|}
\hline
$\Delta x,\Delta t$  & $1/10$ & $1/20$ & $1/40$ & $1/80$  & $1/160$ \tabularnewline\hline

$\Vert \hat{p}_h\Vert_{P}$ & $4.38\times 10^{-2}$ & $3.95\times 10^{-2}$ & $4.20\times 10^{-2}$ & $4.31\times 10^{-2}$ &  $4.33\times 10^{-2}$\tabularnewline

$\Vert \hat{p}_h-p\Vert_{P}$ & $1.80\times 10^{-1}$ & $6.30\times 10^{-2}$ & $1.66\times 10^{-2}$ & $2.78\times 10^{-3}$  & - \tabularnewline

$\Vert \hat{v}_h\Vert_{L^2(0,T)}$  & $1.48\times 10^{-1}$ & $1.33\times 10^{-1}$ & $1.53\times 10^{-1}$ & $1.64\times 10^{-1}$  & $1.67\times 10^{-1}$  \tabularnewline

$\Vert \hat{v}_{h}-v\Vert_{L^2(0,T)}$  & $9.81\times 10^{-2}$ & $6.28 \times 10^{-2}$ & $3.80 \times 10^{-2}$ & $1.11\times 10^{-2}$ & - 
\tabularnewline

$\Vert \hat{y}_h(\cdot\,,T)\Vert_{L^2(0,1)}$  & $1.09\times 10^{-1}$ & $7.67\times 10^{-2}$ & $3.70 \times 10^{-2}$ & $1.11 \times 10^{-2}$ & $1.87\times 10^{-3}$ \tabularnewline

$\Vert \hat{y}_{t,h}(\cdot\,,T)\Vert_{H^{-1}(0,1)}$  & $1.36\times 10^{-1}$ & $8.82 \times 10^{-2}$ & $5.16 \times 10^{-2}$ & $1.76 \times 10^{-2}$ & $2.82\times 10^{-3}$ 
\tabularnewline\hline
\end{tabular}
\caption{$(y_0(x),y_1(x))\equiv (e^{-500(x-0.2)^2},0)$ and $a\equiv 1, b\equiv 1$ - $T=2.2$.}
\label{tab:ex2}
\end{table}

\subsection{Initial data $(y_0,y_1)\in H^1(0,1)\times L^2(0,1)$ and constant speed of propagation}\label{sectionh}

   Let us enhance that our approach, in agreement with the theoretical results, also provides convergent results for irregular initial data. We take a continuous but not differentiable initial state $y_0$ and a piecewise constant initial speed $y_1$: 
   \begin{equation}
y_0(x) \equiv x \, 1_{[0,1/2]}(x) + (1-x) \, 1_{]1/2,1]}(x), \quad y_1(x) \equiv 10 \times 1_{[1/5,1/2]}(x).   \label{irregular_data}
   \end{equation}
The other data are unchanged, except $b$, that is taken equal to zero.

   Observe that these functions remain compatible with the $C^1$ finite element used to approximate~$p$, since $y_0$ and $y_1$ only appear in the right hand side of the variational formulation and~$\pi_{\Delta x} y_0$ and~$\pi_{\Delta x} y_1$ make sense; see~(\ref{pb_phh}). The unique difference is that, once $\hat{p}_h$ and $\hat{v}_h$ are known, $\hat{y}_h$ must be computed from (\ref{post_y1})--(\ref{post_y2}) using a $C^0$ (and not $C^1$) spatial finite element method.
   
   Recall however that these initial data typically generate pathological numerical behavior when the usual dual approach, based on the minimization of~\eqref{I_L2norm}, is used.

\begin{figure}[http]
\begin{center}
\begin{minipage}[t]{6.9cm}
\hspace*{-1.cm}
\includegraphics[scale=0.44]{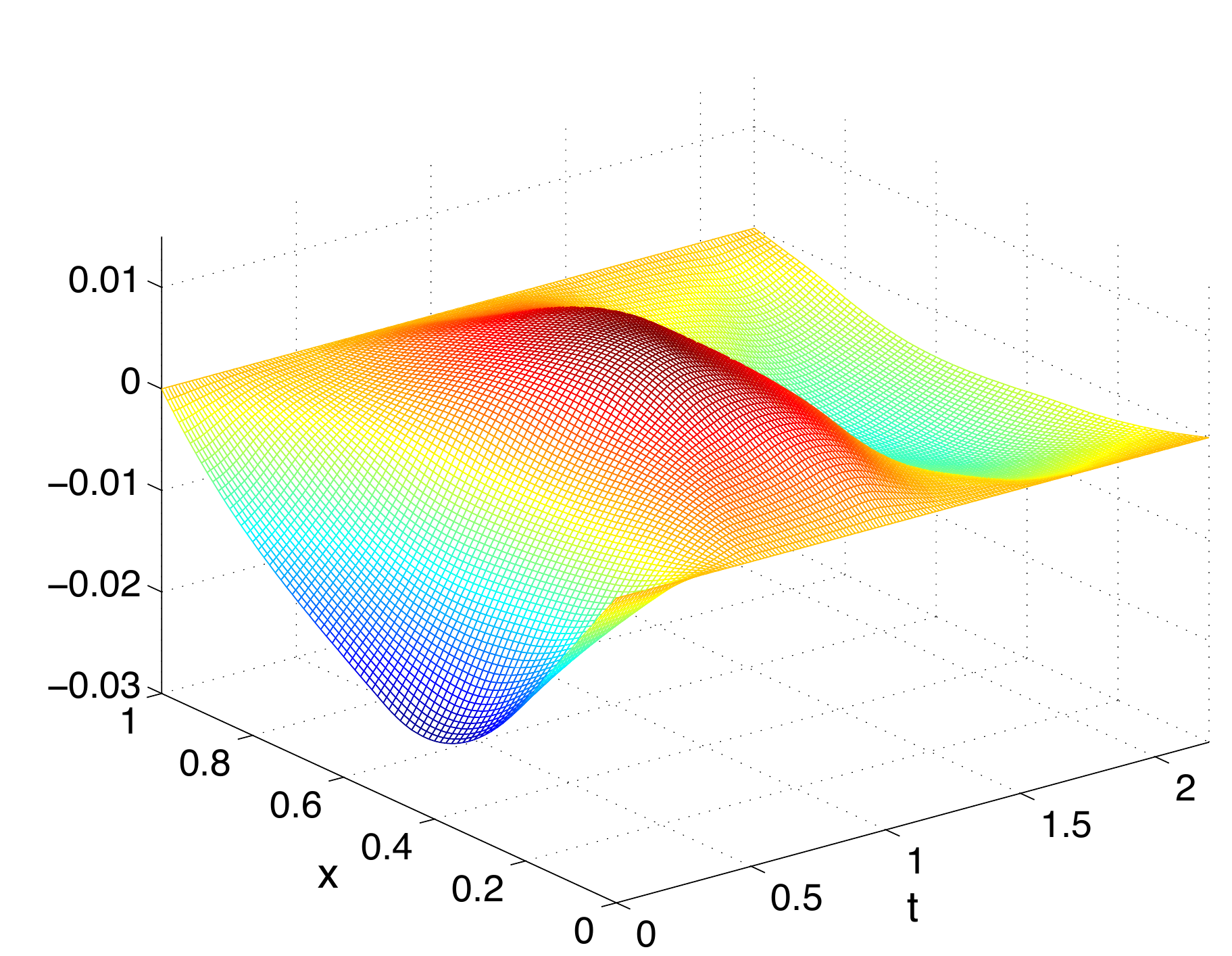}
\end{minipage}
\begin{minipage}[t]{6.9cm}
\hspace*{0.5cm}
\includegraphics[scale=0.38]{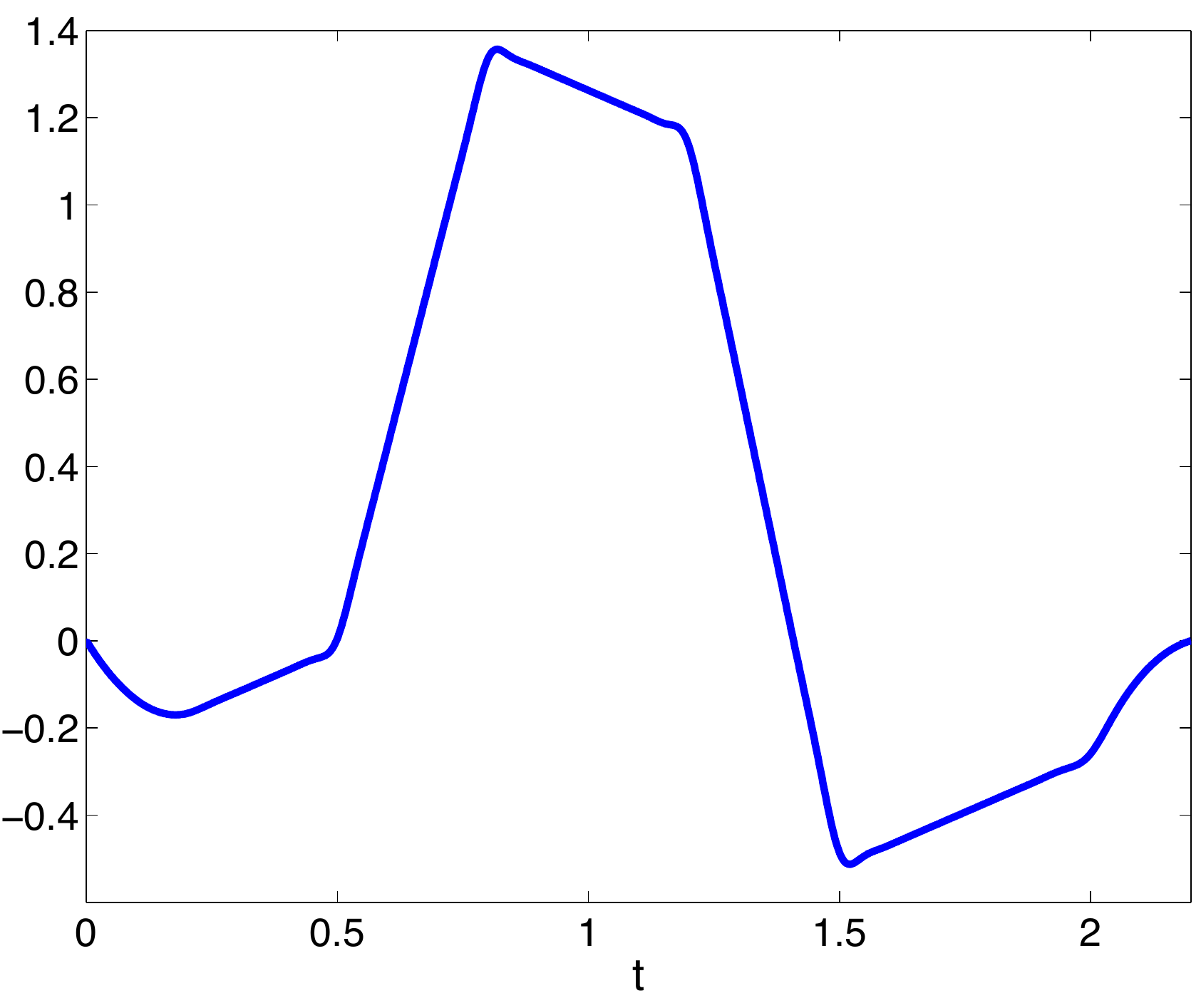}
\end{minipage}
\caption{$(y_0,y_1)$ given by (\ref{irregular_data}) and $a\equiv 1$- The solution $\hat{p}_h$  over $Q_T$ (\textbf{Left}) and the corresponding variable $\hat{v}_h$ on $(0,T)$ (\textbf{Right}) - $h=(1/80,1/80)$.}\label{fig:p_controle_ex3}
\end{center}
\end{figure}

\begin{figure}[http!]
\begin{center}
\includegraphics[scale=0.5]{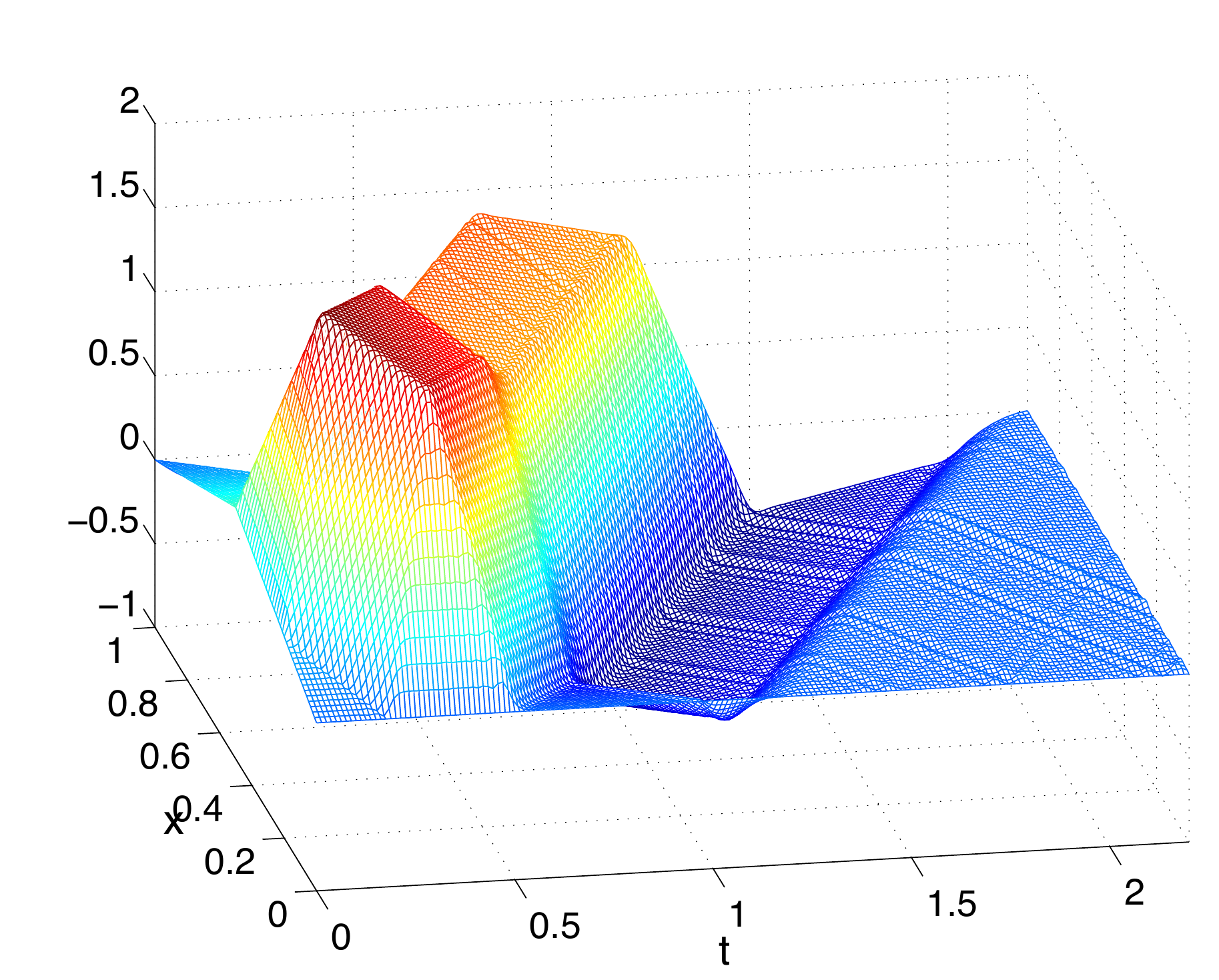}
\caption{$(y_0,y_1)$ given by (\ref{irregular_data}) and $a\equiv 1$ -The solution $\hat{y}_h$  over $Q_T$ - $h=(1/80,1/80)$.}\label{fig:y_ex3}
\end{center}
\end{figure}

   Some numerical results are given in~Table~\ref{tab:ex3} and~Figures~\ref{fig:p_controle_ex3} and~\ref{fig:y_ex3}. As before, we observe the convergence of the variable $\hat{p}_h$ and therefore $\hat{v}_{h}$ and $\hat{y}_h$ as $h \to 0$. We see that $\Vert \hat{p}_h-p\Vert_P = \mathcal{O}(h^{1.48})$ and $\Vert \hat{v}_h-v\Vert_{L^2(0,1)} = \mathcal{O}(h^{1.23})$. In particular, we do not observe oscillations for the control or the functions $\hat{p}_h$ and $\hat{p}_{h,t}$ at the initial time.

\

\begin{table}[http]
\centering
\begin{tabular}{|c|ccccc|}
\hline
$\Delta x,\Delta t$  & $1/10$ & $1/20$ & $1/40$ & $1/80$  & $1/160$ \tabularnewline\hline

$\Vert \hat{p}_h\Vert_P$ & $3.16\times 10^{-1}$ & $2.89\times 10^{-2}$ & $2.73\times 10^{-2}$ & $2.65\times10^{-2}$ &  $2.61\times 10^{-1}$\tabularnewline

$\Vert \hat{p}_h-p\Vert_P$ & $1.12\times 10^{-1}$ & $4.62\times 10^{-2}$ & $1.70\times 10^{-2}$ & $5.12\times 10^{-3}$  & - \tabularnewline

$\Vert \hat{v}_h\Vert_{L^2(0,T)}$  & $1.23$ & $1.11$ & $1.05$ & $1.02$  & $1.004$  \tabularnewline

$\Vert \hat{v}_{h}-v\Vert_{L^2(0,T)}$  & $2.52\times 10^{-1}$ & $1.25\times 10^{-1}$ & $5.57\times 10^{-2}$ & $1.90\times 10^{-2}$ & - 
\tabularnewline

$\Vert \hat{y}_h(\cdot\,,T)\Vert_{L^2(0,1)}$  & $1.09\times 10^{-1}$ & $5.40\times 10^{-2}$ & $2.20\times 10^{-2}$ & $1.09\times 10^{-2}$ & $6.20\times 10^{-3}$ \tabularnewline

$\Vert \hat{y}_{t,h}(\cdot\,,T)\Vert_{H^{-1}(0,1)}$  & $7.25\times 10^{-2}$ & $4.62\times 10^{-2}$ & $2.85\times 10^{-2}$ & $5.12\times 10^{-3}$ & $6.75\times 10^{-3}$ 
\tabularnewline\hline
\end{tabular}
\caption{$(y_0,y_1)$ given by (\ref{irregular_data}) and $a\equiv 1$ - $T=2.2$.}
\label{tab:ex3}
\end{table}

\subsection{Discontinuous initial data $y_0$ and constant speed of propagation}

\begin{figure}[http]
\begin{center}
\begin{minipage}[t]{6.9cm}
\hspace*{-1.cm}
\includegraphics[scale=0.44]{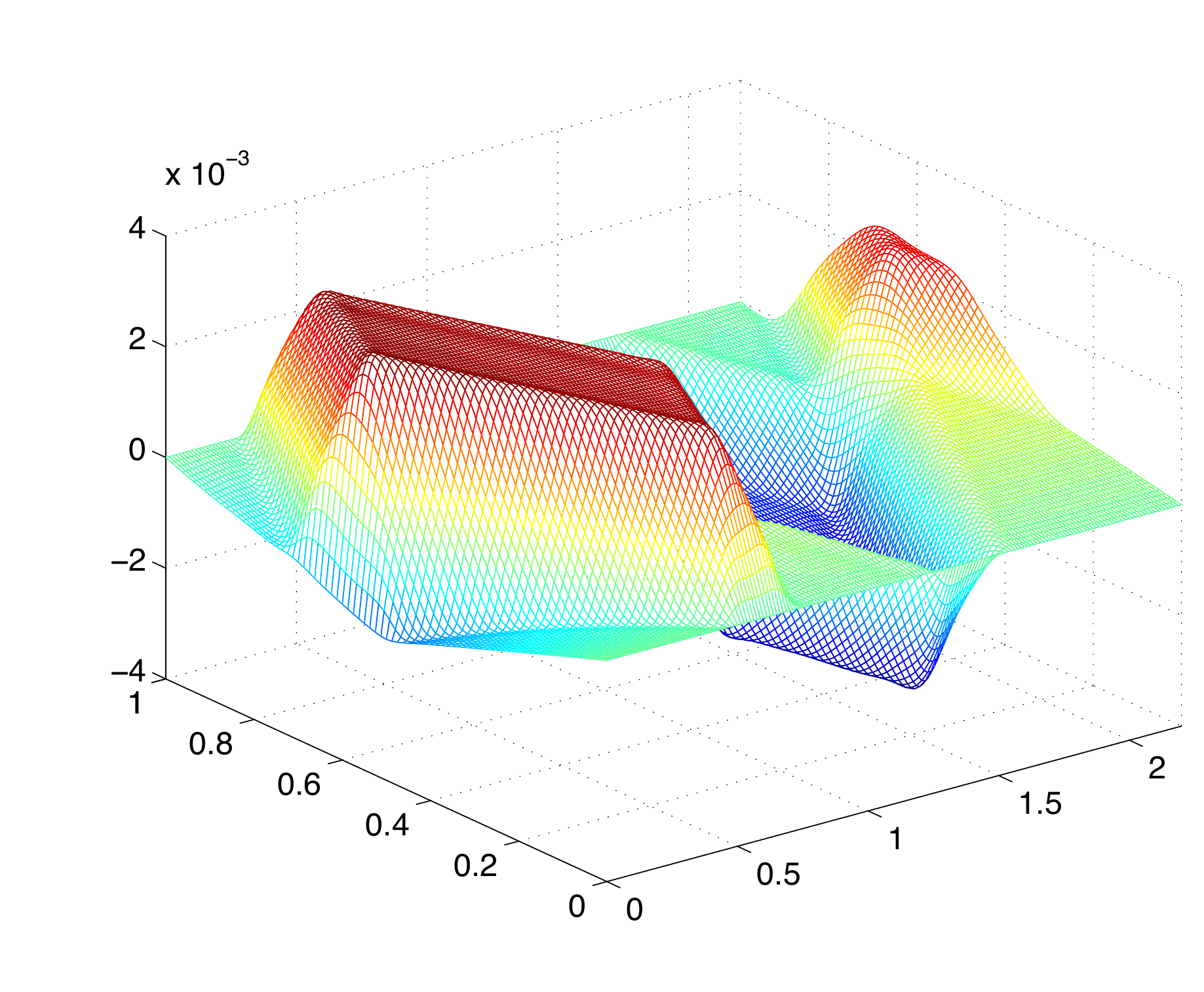}
\end{minipage}
\begin{minipage}[t]{6.9cm}
\hspace*{0.5cm}
\includegraphics[scale=0.38]{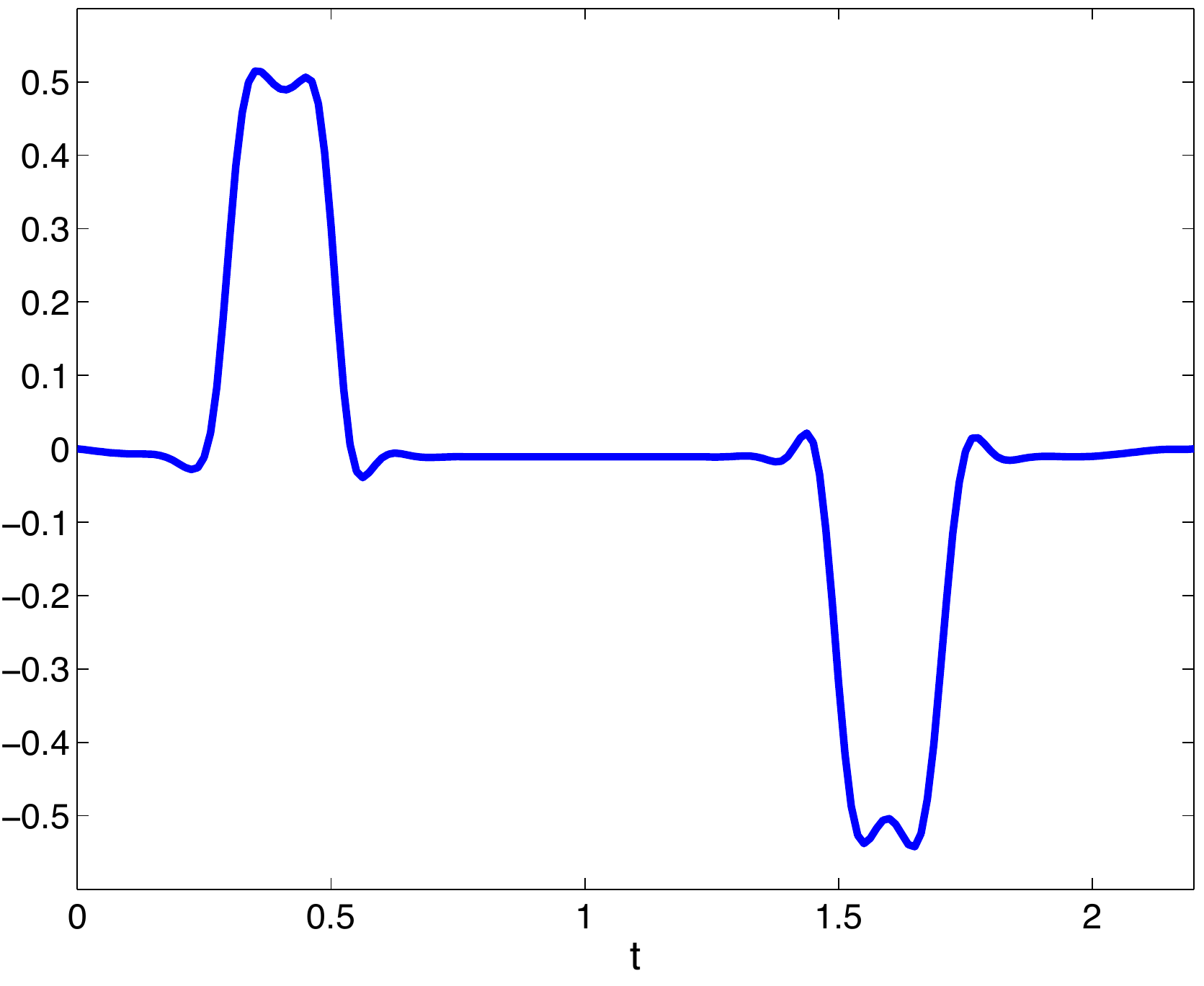}
\end{minipage}
\caption{$(y_0,y_1)$ given by (\ref{very_irregular_data}) and $a\equiv 1$- The solution $\hat{p}_h$  over $Q_T$ (\textbf{Left}) and the corresponding variable $\hat{v}_h$ on $(0,T)$ (\textbf{Right}) - $h=(1/80,1/80)$.}\label{fig:p_controle_ex3b}
\end{center}
\end{figure}

The method also provides convergent results for data $y_0$ only in $L^2(0,1)$. We consider the following initial condition:
 \begin{equation}
y_0(x) \equiv 1_{[0.5,0.7]}(x), \quad y_1(x) \equiv 0.   \label{very_irregular_data}
   \end{equation}
The other data are unchanged with respect to Section \ref{sectionh}. This leads to pathological numerical behavior when other frequently used dual methods are employed (we refer to \cite{munch05}).  Some numerical results are given in~Table~\ref{tab:ex3b} and~Figure~\ref{fig:p_controle_ex3b}. Once again, the convergence of the variable $\hat{p}_h$ and therefore $\hat{v}_{h}$ and $\hat{y}_h$ as $h \to 0$ is observed. 

\begin{table}[http]
\centering
\begin{tabular}{|c|ccccc|}
\hline
$\Delta x,\Delta t$  & $1/10$ & $1/20$ & $1/40$ & $1/80$  & $1/160$ \tabularnewline\hline

$\Vert \hat{p}_h\Vert_P$ & $1.01\times 10^{-1}$ & $1.00\times 10^{-1}$ & $9.71\times 10^{-2}$ & $9.53\times10^{-2}$ &  $9.47\times 10^{-2}$\tabularnewline

$\Vert \hat{v}_h\Vert_{L^2(0,T)}$  & $3.42\times 10^{-1}$ & $3.27\times 10^{-1}$ & $3.19\times 10^{-1}$ & $3.14\times 10^{-1}$  & $3.14\times 10^{-1}$  \tabularnewline

$\Vert \hat{y}_h(\cdot\,,T)\Vert_{L^2(0,1)}$  & $1.24\times 10^{-1}$ & $9.27\times 10^{-2}$ & $7.26\times 10^{-2}$ & $5.88\times 10^{-2}$ & $3.12\times 10^{-2}$ \tabularnewline

$\Vert \hat{y}_{t,h}(\cdot\,,T)\Vert_{H^{-1}(0,1)}$  & $1.55\times 10^{-1}$ & $1.16\times 10^{-1}$ & $1.06\times 10^{-1}$ & $7.13\times 10^{-2}$ & $6.02\times 10^{-2}$ 
\tabularnewline\hline
\end{tabular}
\caption{$(y_0,y_1)$ given by (\ref{very_irregular_data}) and $a\equiv 1$ - $T=2.2$.}
\label{tab:ex3b}
\end{table}

\subsection{Non constant smooth speed of propagation}

   {\color{black} Finally,  let us consider a non-constant function $a = a(x)$ (we refer to~\cite{Glo02} for the dual approach in this case). In order to illustrate the robustness of our method, we will take a coefficient $a \in C^1([0,1])$ } with
   \begin{equation}
a(x)=
\left\{
\begin{aligned}
&  1 &  x\in [0,0.45] \\
& \in [1.,5.] \quad (a^{\prime}(x)>0),  & x\in (0.45,0.55) \\
&  5 & x\in [0.55,1]
\end{aligned}
\right.
\label{a_nconstant}
   \end{equation} 
so that condition (\ref{eq:Tlargebis}) is equivalent to $T>2(1+1/20)\sqrt{5}\approx 4.69$ (taking again $x_0=-1/20$). In order to reduce the computational cost, we take as before $T=2.2$ and we still observe that the constant $C_{0h}$ in~\eqref{43p} is uniformly bounded. 
 
\begin{figure}[http]
\begin{center}
\begin{minipage}[t]{6.9cm}
\hspace*{-1.cm}
\includegraphics[scale=0.44]{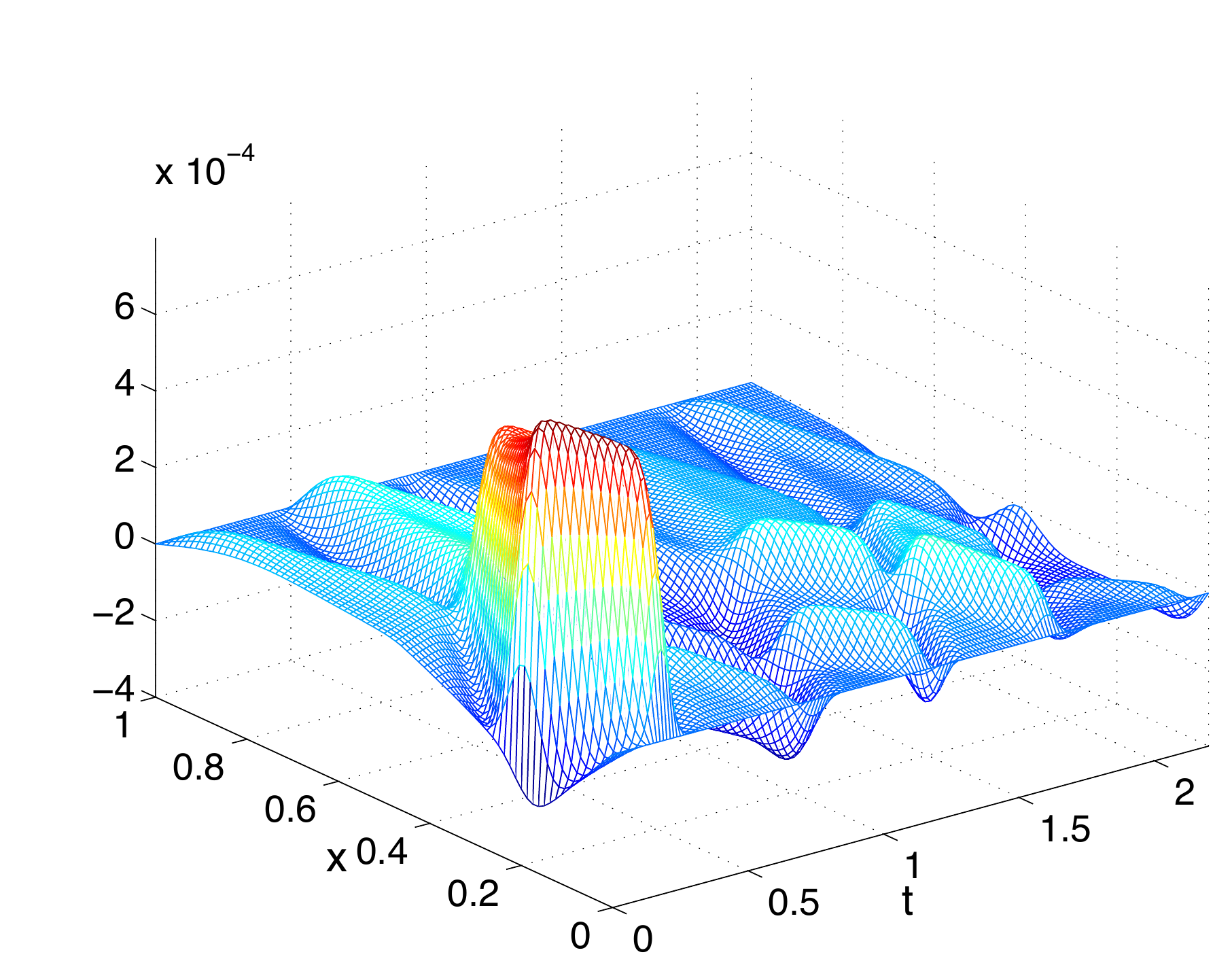}
\end{minipage}
\begin{minipage}[t]{6.9cm}
\hspace*{0.5cm}
\includegraphics[scale=0.38]{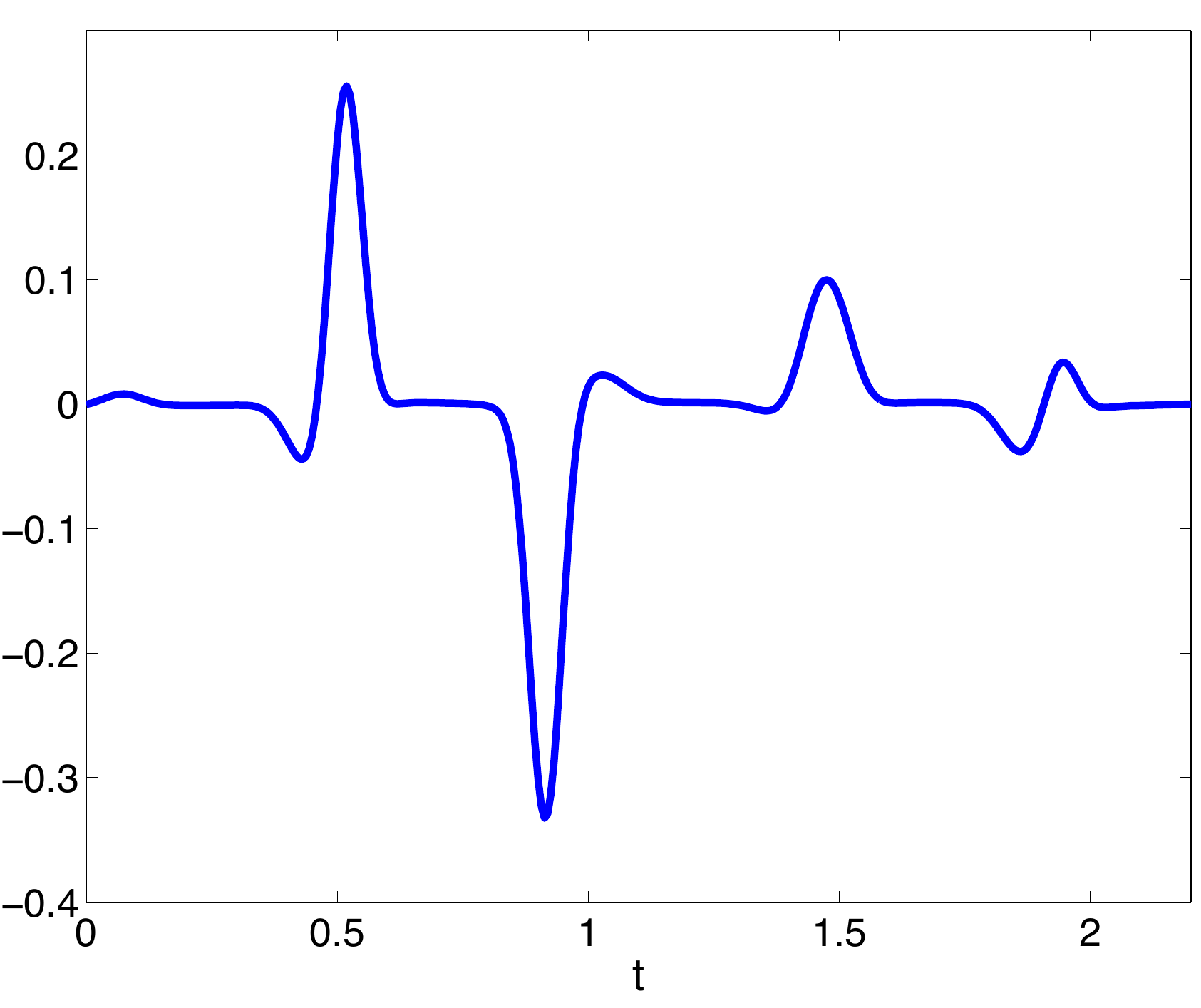}
\end{minipage}
\caption{$(y_0(x),y_1(x))\equiv (e^{-500(x-0.2)^2},0)$ and $a$ given by (\ref{a_nconstant}) - The solution $\hat{p}_h$  over $Q_T$ (\textbf{Left}) and the corresponding variable $\hat{v}_h$ on $(0,T)$ (\textbf{Right}) - $h=(1/80,1/80)$.}\label{fig:p_controle_ex4}
\end{center}
\end{figure}

\begin{figure}[http!]
\begin{center}
\includegraphics[scale=0.5]{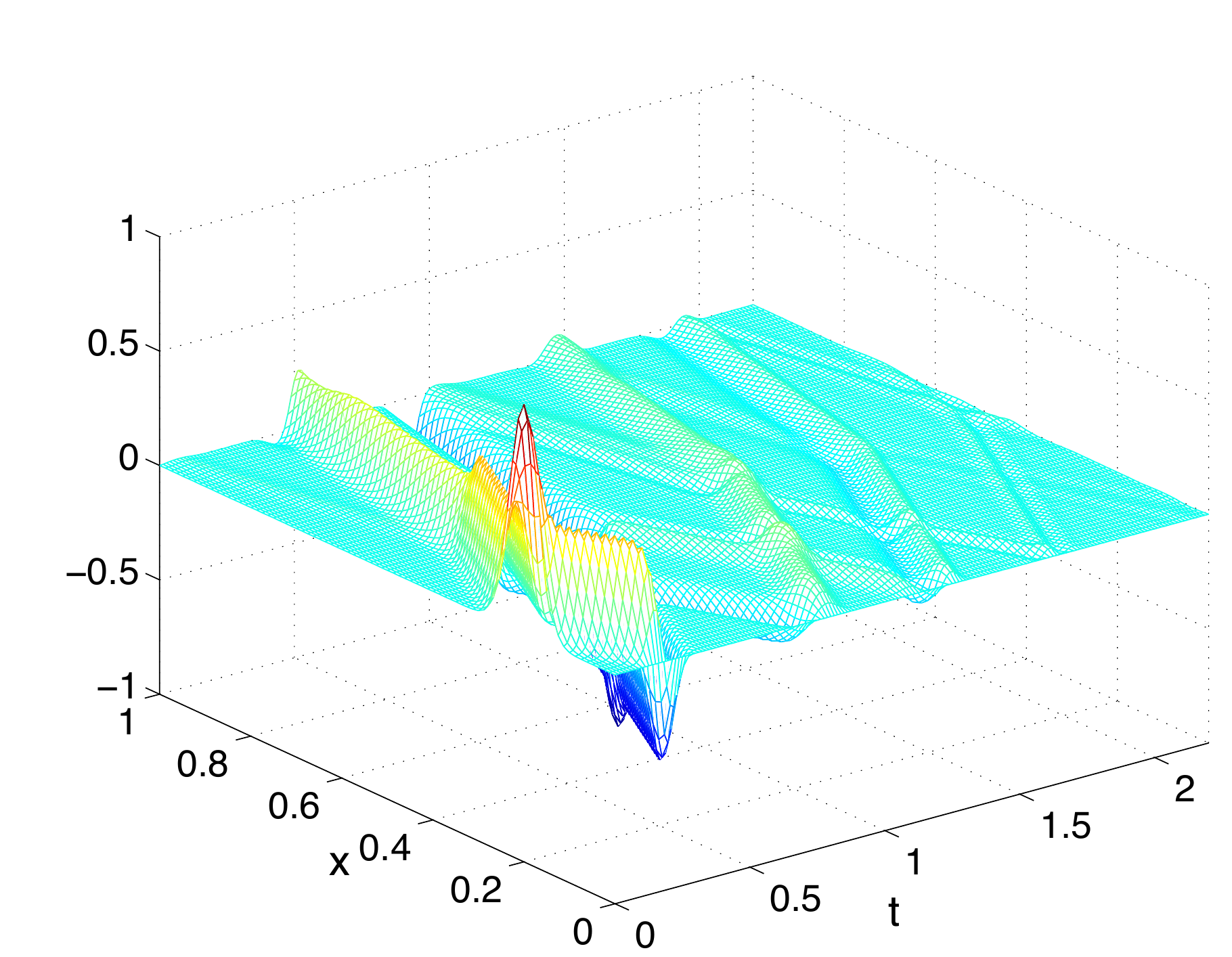}
\caption{$y_0(x)\equiv e^{-500(x-0.2)^2}$ and $a$ given by (\ref{a_nconstant}) -The solution $\hat{y}_h$  over $Q_T$ - $h=(1/80,1/80)$.}\label{fig:y_ex4}
\end{center}
\end{figure}

We take again $(y_0(x),y_1(x)) \equiv (e^{-500(x-0.2)^2},0)$ and $b \equiv 0$. Table~\ref{tab:ex4} illustrates the convergence of the approximations with respect to $h$. Figures~\ref{fig:p_controle_ex4} and~\ref{fig:y_ex4} depict for $h=(1/80,1/80)$ the functions $\hat{p}_h$, $\hat{v}_h$ and $\hat{y}_h$. In particular, in~Figure~\ref{fig:y_ex4}, we can observe the diffraction of the wave when crossing the transitional zone $(0.45,0.55)$.

\begin{table}[http]
\centering
\begin{tabular}{|c|ccccc|}
\hline
$\Delta x,\Delta t$  & $1/10$ & $1/20$ & $1/40$ & $1/80$  & $1/160$ \tabularnewline\hline

$\Vert \hat{p}_h\Vert_{P_h}$ & $3.87\times 10^{-2}$ & $3.44\times 10^{-2}$ & $3.75\times10^{-2}$ & $3.85\times10^{-2}$ &  $3.86\times 10^{-2}$\tabularnewline

$\Vert \hat{p}_h-p\Vert_{P_h}$ & $1.25\times 10^{-1}$ & $5.75\times 10^{-2}$ & $2.64\times10^{-2}$ & $1.01\times10^{-2}$  & - \tabularnewline

$\Vert \hat{v}_h\Vert_{L^2(0,T)}$  & $7.74\times 10^{-2}$ & $6.53\times 10^{-2}$ & $9.16\times10^{-2}$ & $1.01\times10^{-1}$  & $1.03\times 10^{-1}$  \tabularnewline

$\Vert \hat{v}_{h}-v\Vert_{L^2(0,T)}$  & $5.07\times 10^{-1}$ & $4.17\times10^{-2}$ & $2.03\times10^{-2}$ & $4.86\times10^{-3}$ & - 
\tabularnewline

$\Vert \hat{y}_h(\cdot\,,T)\Vert_{L^2(0,1)}$  & $1.09\times 10^{-1}$ & $7.89\times10^{-2}$ & $1.81\times10^{-2}$ & $1.16\times10^{-2}$ & $1.71\times 10^{-3}$ \tabularnewline

$\Vert \hat{y}_{t,h}(\cdot\,,T)\Vert_{H^{-1}(0,1)}$  & $1.01\times10^{-1}$ & $8.39\times10^{-2}$ & $4.81\times10^{-2}$ & $7.52\times10^{-3}$ & $1.55\times 10^{-3}$ 
\tabularnewline\hline
\end{tabular}
\caption{$(y_0(x),y_1(x))\equiv (e^{-500(x-0.2)^2,0)}$ and $a$ given by (\ref{a_nconstant}) - $T=2.2$.}
\label{tab:ex4}
\end{table}

\section{Further comments and concluding remarks}
\label{sec:end_remarks}

   Let us begin this section with some general considerations on the use of Carleman weights that serve to justify our approach:
   
\begin{enumerate}
%
\item The search of a control minimizing $J$ in~\eqref{P-FI}, where $y$ is involved, is very appropriate from the numerical viewpoint. As shown in Section \ref{sec:var}, the explicit occurrence of the state variable $y$ leads to an elliptic problem in~$Q_T$, that is easy to analyze and solve (at this level, the particular choice of the weight is less important). This approach does not require the discretization of the wave operator, as for usual dual approachs; therefore, it does not generate any spurious oscillations and leads to numerical well-posedness. This is an important feature of the approach.   
\item The Carleman weights provide regularity of the solution to~\eqref{eq:varp} and therefore allows to derive estimates of the errors $\Vert p-p_h\Vert_P$ in term of $h=(\Delta x,\Delta t)$. This will be detailed in a forthcoming work.
\item The process can be viewed as a first step for the numerical controllability of semi-linear problems:
  if we just apply a fixed-point argument, we will find at each iterate a linear equation with non-regular coefficients depending on $x$ and~$t$ for which the present approach is adequate.
\item In our numerical experiments we have not found any essential difference for small or large $s$ or $\lambda$: this is in full agreement with Remark \ref{re:about_s} and Section \ref{carlemanconstant}.
\item Furthermore, we mention that the approach has been considered by analogy with the similar analysis in references~\cite{EFC-AM-11} to~\cite{EFC-AM-long}, dealing with heat equations.
\end{enumerate}


\subsection{Primal versus dual approach~(I): analogies}

   The solution to the variational formulation (\ref{eq:varp}) is also the unique minimizer of the functional~$I$, with 
   \begin{equation}\label{newI}
   \begin{array}{l}
\dis I(p):= \frac{1}{2}\jjntQT \rho^{-2}\vert Lp\vert^2\,dx\,dt + \frac{1}{2}\int_0^T \rho_0^{-2} a(1)^2 \vert p_x(1,\cdot)\vert^2\, dt \\
\dis \phantom{I(p) } \quad -\int_0^1 y_0(x)\,p_t(x,0)\,dx + \langle y_1,p(\cdot,0)\rangle_{H^{-1},H_0^1}.
   \end{array}
   \end{equation}
This is similar to the conjugate functional $J^{\star}$ in~(\ref{P-FIstar}). Actually, we notice that 
  $J^{\star}(\mu,\phi_0,\phi_1)= I(-\phi)$ for all $(\mu,\phi_0,\phi_1)\in L^2(Q_T)\times \boldsymbol{H}$. 
   
   Therefore, the extremal problems (\ref{P-FIstar}) and (\ref{newI}) are connected to each other having (\ref{P-FI}) as starting point. The problem 
   (\ref{eq:varp}), deduced from the primal approach belongs to the framework of elliptic variational problems in two dimensions and is well tailored for a resolution with finite elements.  The dual problem (\ref{P-FIstar}) is of hyperbolic nature: the time variable is kept explicitly and time integration is required. 
   
   Note that we may also derive the optimality conditions for $J^{\star}$ (as we did in Section \ref{sec:var} for $J$): this leads, at least formally, to the problem (\ref{eq:varp}).
   
   We also mention \cite{pedregal} where a (different) variational approach is introduced. 


\subsection{Primal versus dual approaches~(II): discrete properties}

   The variational approach used here leads to satisfactory convergence results, in particular the strong convergence of the approximate controls $\hat{v}_h$ towards a null control of the wave equation. This relies in a fundamental way on the fact that we work in a subspace $P_h$ of~$P$. Indeed, this allows to write directly the Carleman estimate in $P_h$ and get that the function $I$ (given by (\ref{newI})) is uniformly coercive with respect to the discretization parameter $h$.
   
   On the other hand, notice that no wave equation has to be solved in order to compute the approximations $\hat{v}_h$. For each $h$, once $\hat{v}_h$ is known, we must solve the wave equation, in a post-treatment process, to compute the corresponding state $\hat{y}_h$ (recall that, actually, this may be avoided by using directly the optimality condition $y=-\rho^{-2}Lp$).

   This is in contrast with the dual approach. Indeed, the minimization of $J^{\star}$ by an iterative process requires the resolution of wave equations, through a decoupled space and time discretization. As recalled in the introduction, this may lea to numerical pathologies (the occurrence of spurious high frequency solutions) and, therefore, needs some specific numerical approximations and techniques. We mention the work~\cite{Bau11-bis}, where the authors prove, in a close context and within a dual approach, a weaker uniform semi-discrete Carleman estimate with an additional term in the right hand side, necessary to absorb these possibly spurious high frequencies  (see~\cite{Bau11-bis}, Theorem~2.3). 

   Notice that the computed $\hat{v}_h$ are not \textit{a priori} null controls for discrete systems (associated to the wave equation (\ref{eq:wave})), but simply approximations of the control $v$ furnished by the solution to~\eqref{P-FI}. {\color{black} If one wants to go further in the comparison, it can be said that the primal approach aims to first compute the control for (\ref{eq:wave}) and then approximate it, while the dual classical method aims first to discretize (\ref{eq:wave}) and then control the corresponding finite dimensional system.}

   Let us also observe that the (primal) approach in this paper is relatively easy to implement. In practice, the resolution is reduced to solve a linear system, with a banded sparse, symmetric and definite positive matrix, for which efficient direct $LU$ type solvers are known and available. Furthermore, we may want to adapt (and refine locally) the mesh of $Q_T$ in order to improve convergence and such adaptation is much simpler than in the dual approach, where $t$ is "conserved" as a time variable.
  For additional considerations, see also~\cite{Deh-Leb} and~\cite{Leb-Nod}.
   

\subsection{Mixed formulation and $C^0$-approximation}

   The approach can be extended to the higher dimensional case of the wave equation in a bounded set $\Omega \subset \mathbb{R}^N$, with~$N \geq 2$. However, the use of $C^1$-finite element is a bit more involved. Arguing as in~\cite{EFC-AM-11}, we may avoid this difficulty by introducing a mixed formulation equivalent to (\ref{eq:varp}).
   
   The idea is to keep explicit the variable $y$ in the formulation and to introduce a Lagrange multiplier, associated to the constraint $\rho^2 y+Lp=0$ (see~(\ref{eq:yvp})). We obtain the following mixed formulation: find $(y,p,\lambda) \in Z\times P \times Z$ such that 
   \begin{equation}
\label{mixed}
\left\{
\begin{aligned}
&\jjntQT \rho^{2} y\,\overline{y}\,dx\,dt + \int_0^T \rho_0^{-2} a(1)^2 p_x(1,t) \overline{p}_x(1,t)\,dt
+ \jjntQT \lambda (\rho^2 \overline{y} + L\overline{p})\,dx\,dt \\
& \qquad = \int_0^1 y_0(x)\, \overline{p}_t(x,0)\,dx 
- \langle y_1,\overline{p}(\cdot,0)\rangle_{H^{-1}, H_0^1} \quad \forall (\overline{y},\overline{p})\in Z\times P,\\
& \jjntQT \overline{\lambda}(\rho^2 y+Lp)\,dx\,dt=0 \quad \forall \overline{\lambda}\in Z,
\end{aligned}
\right.
   \end{equation}
where
   $$
Z = L^2(\rho^2; Q_T):= \{\, z\in L^1_{loc}(Q_T) :  \jjntQT  \rho^2 \vert z\vert^2\, dx\,dt<+\infty \,\}.
   $$
Taking advantage of the global estimate (\ref{eq:CarlT}), we may show, through an appropriate inf-sup condition, that (\ref{mixed}) is well-posed in $Z\times P\times Z$. Moreover, the approximation of this formulation may be addressed using $C^0$-finite element, which is very convenient. The approximation is non-conformal. More precisely, the variable $p$ is now sought in a space $R_h$ of $C^0$-functions that is not included in $P$.

   At the discrete level, (\ref{mixed}) reduces the controllability problem to the inversion of a square, banded and symmetric matrix. Moreover, as before, no wave equation has to be solved, whence the numerical pathology described above is not expected. However, since the underlying approximation is not conformal (this is the price to pay to avoid $C^1$ finite elements), a careful (and \textit{a priori} not straightforward) choice for $R_h$ has to be done in order to guarantee a uniform discrete inf-sup condition. The analysis of this point, as well as the use of stabilized finite elements, will be detailed in a future work. 

\


\subsection{Extensions}

   The approach presented here can be extended and adapted to other equations and systems. What is needed is, essentially, an appropriate Carleman estimate. 
   
   In particular, we can adapt the previous ideas and results to the inner controllability case, i.e.~the null controlability of the wave equation with distributed controls acting on a (small) sub-domain $\omega$ of $(0,1)$. Furthermore, using finite element tools, we can also get results in the case where the sub-domain $\omega$ varies in time, that is non-cylindrical control domains $q_T$ of the form
   $$
q_T=\{\,(x,t)\in Q_T: g_1(t)<x<g_2(t), \quad t\in (0,T)\,\},
   $$
where $g_1$ and $g_2$ are smooths functions on $[0,T]$, with $0\leq g_1 < g_2 \leq 1$. This opens the possibility to optimize numerically the domain $q_T$, as was done in a cylindrical situations in \cite{munch08} (see also \cite{periago}).

   Let us finally mention that many non-linear situations can be considered through a suitable linearization and iterative process. We refer to~\cite{EFC-AM-NL, EFC-AM-long} for some ideas in a similar parabolic situation. 


\appendix
\section{Appendix: On the proof of Theorem 2.1}

   We first prove a global Carleman estimate for functions $w$ satisfying vanishing initial and final conditions. In what follows, $L$ stands for the operator given in~\eqref{eq:L} with $b \equiv 0$. It is easy to check that, if the estimate \eqref{eq:CarlT} holds in this particular case, then the same estimate holds for any potential $b \in L^\infty((0,1) \times (-T,T))$.

\begin{theorem}\label{th:CarlemanV}
    With the notation of Section~\ref{sec:var}, let $x_0 < 0$ be a fixed point, let $\phi$ and $\varphi$ be the weight functions defined by \eqref{eq:phi}--\eqref{eq:varphi} and let $a \in \mathcal{A}(x_0, a_0)$ with $a_0 > 0$. Then there exist positive constants $s_0$ and $M$, only depending on $x_0$, $a_0$, {\color{black} $\|a\|_{C^3([0,1])}$ } and~$T$ such that, for all $s > s_0$, one has:
   \begin{equation}\label{eq:CarlTbis}
\begin{array}{c}
\dis s \int_{-T}^T \int_0^1 e^{2s\varphi} \left( |v_t|^2 + |v_x|^2 \right) \,dx\,dt +  s^3 \int_{-T}^T  \int_0^1 e^{2s\varphi} |v|^2 \,dx\,dt \\
\noalign{\smallskip}
\dis \leq M \int_{-T}^T \int_0^1 e^{2s\varphi} |L v |^2 \,dx\,dt  + M s \int_{-T}^{T} e^{2s\varphi} |v_x(1,t)|^2 \, dt
\end{array}
   \end{equation}
for any $v \in L^2(-T, T; H_0^1(0,1))$ satisfying $L v \in L^2((0, 1) \times (-T, T))$, $v_x(1, \cdot) \in L^2(-T, T)$ and
   \[
v(\cdot\,,\pm T) = v_t(\cdot\,,\pm T) = 0.
   \]
\end{theorem}

   The proof of this result follows step-by-step the proof of Theorem 2.1 in \cite{Bau11}. However, since the argument provides conditions on the set of admissible $a$ and, to our knowledge, these conditions have not been stated in this form before, we provide here the detailed proof.
   
\

\noindent
\textsc{Proof:} Let us introduce $w = e^{s \varphi} v$ and let us set
   $$
P w:= e^{s \varphi} L(e^{-s \varphi} w) = e^{s \varphi} \left((e^{-s \varphi} w)_{tt} - (a (e^{-s \varphi} w)_x)_x\right).
   $$
After some computations, we find that $Pw = P_1 w + P_2 w + Rw$, with
   \begin{align*}
  & P_1 w = w_{tt} - (a w_x)_x + s^2 \lambda^2 \varphi^2 w \left( |\psi_t|^2 - a |\psi_x|^2 \right) \\
  & P_2 w = (\alpha - 1) s \lambda \varphi w \left( \psi_{tt} - (a \psi_x)_x \right) - s\lambda^2 \varphi w \left( |\psi_t|^2 - a |\psi_x|^2 \right) -2s\lambda\varphi \left( \psi_t w_t - a \psi_x w_x \right) \\
  & Rw = - \alpha s \lambda \varphi w \left( \psi_{tt} - (a \psi_x)_x \right),
   \end{align*}
where the parameter $\alpha$ will be chosen below.

   Recall that
   $$
\psi(x, t) \equiv |x-x_0|^2 - \beta t^2 + M_0, \quad \varphi(x, t) \equiv e^{\lambda \psi(x, t)}
   $$
and
   $$
\psi(x, t) \geq 1 \quad \forall (x, t) \in (0, 1) \times (-T, T).
   $$
   In this proof, we will denote by $M$ a generic positive constant that can depend on~$x_0$, $a_0$, {\color{black} $\Vert a\Vert_{C^3([0,1])}$} and~$T$.
 
   As in the constant case $a \equiv 1$, the first part of the proof is devoted to estimate from below the integral
   \begin{equation}
\label{eq:P1P2}
I = \int_{-T}^{T} \int_\Omega (P_1 w) \, (P_2 w) \,dx\,dt = \sum_{i, j = 1}^3 I_{ij}.
   \end{equation}
   
  By integrating by parts in time and/or space, we can compute the integrals $I_{ij}$ in~\eqref{eq:P1P2}. We obtain:
\begin{eqnarray*}
  I_{11} & = & (\alpha - 1)s \lambda\int_{-T}^{T} \int_0^1 w_{tt} \, \varphi w( \psi_{tt} - (a\psi_x)_x )\,dx\,dt \\
  &=&(1-\alpha )s\lambda \int_{-T}^{T} \int_0^1 \varphi  | w_{t} | ^2 (\psi_{tt} - (a\psi_x)_x ) \,dx\,dt\\
  &&- \dfrac {(1-\alpha )}2 s\lambda^2 \int_{-T}^{T} \int_0^1 \varphi | w | ^2 \psi_{tt}(\psi_{tt} - (a\psi_x)_x ) \,dx\,dt\\
  &&- \dfrac {(1-\alpha )}2 s\lambda^3 \int_{-T}^{T} \int_0^1 \varphi | w | ^2  | \psi_t |^2(\psi_{tt} - (a\psi_x)_x ) \,dx\,dt,
\end{eqnarray*}
\begin{eqnarray*}
  I_{12}&=&-s\lambda^2 \int_{-T}^{T} \int_0^1 w_{tt} \, \varphi w(| \psi_t|^2-a|\psi_x|^2)\,dx\,dt\\
  &=&s \lambda ^2 \int_{-T}^{T} \int_0^1 \varphi |w_t|^2(| \psi_t|^2-a|\psi_x|^2)\,dx\,dt
  -s\lambda^2 \int_{-T}^{T} \int_0^1 \varphi  |w|^2 |\psi_{tt}|^2\,dx\,dt\\
  &&-\dfrac {3s\lambda^3}2 \int_{-T}^{T} \int_0^1 \varphi |w|^2 |\psi_t|^2 \psi_{tt}\,dx\,dt
  +\dfrac {s\lambda^3}2 \int_{-T}^{T} \int_0^1 \varphi |w|^2 a |\psi_x|^2 \psi_{tt}\,dx\,dt\\
  &&-\dfrac {s\lambda^4}2 \int_{-T}^{T} \int_0^1 \varphi  |w|^2 | \psi_t|^2 (| \psi_t|^2-a|\psi_x|^2)\,dx\,dt
\end{eqnarray*}
and
\begin{eqnarray*}
  I_{13}&=& -2s\lambda \int_{-T}^{T} \int_0^1 w_{tt} \, \varphi (\psi_t w_t - a \psi_x w_x))\,dx\,dt\\
  &=& s\lambda \int_{-T}^{T} \int_0^1 \varphi | w_t | ^2 \psi_{tt} \,dx\,dt
  +s\lambda^2 \int_{-T}^{T} \int_0^1 \varphi | w_t | ^2 | \psi_t |^2 \,dx\,dt\\
  &&+ s\lambda \int_{-T}^{T} \int_0^1 \varphi | w_t | ^2 (a \psi_x)_x  \,dx\,dt
  + s\lambda^2 \int_{-T}^{T} \int_0^1 \varphi | w_t | ^2  a |\psi_x|^2 \,dx\,dt\\
  &&-2s\lambda^2 \int_{-T}^{T} \int_0^1 \varphi a \psi_x \psi_t \, w_x w_t \,dx\,dt.
\end{eqnarray*}

   Also,
\begin{eqnarray*}
  I_{21}&=&(1-\alpha )s \lambda\int_{-T}^{T} \int_0^1 (a w_x)_x \, \varphi w(\psi_{tt} - (a\psi_x)_x )\,dx\,dt\\
  &=& -(1-\alpha ) s \lambda \int_{-T}^{T} \int_0^1 \varphi a |w_x|^2 (\psi_{tt} -(a \psi_x)_x ) \,dx\,dt\\
  &&+ \dfrac {(1-\alpha )}2 s \lambda^2 \int_{-T}^{T} \int_0^1 \varphi |w|^2 (a \psi_x)_x (\psi_{tt} - (a\psi_x)_x )\,dx\,dt\\
  &&+ \dfrac {(1-\alpha )}2 s \lambda^3 \int_{-T}^{T} \int_0^1 \varphi a |w|^2  |\psi_x|^2 (\psi_{tt} -(a \psi_x)_x )\,dx\,dt \\
  & & - (1- \alpha) s \lambda^2 \int_{-T}^{T} \int_0^1 \varphi |w|^2 a\psi_x (a \psi_x)_{xx}  \,dx\,dt \\
  & & - \dfrac {(1-\alpha )}2 s \lambda \int_{-T}^{T} \int_0^1 \varphi |w|^2 \left(a_x (a \psi_x)_{xx} + a (a \psi_x)_{xxx} \right) \,dx\,dt,
\end{eqnarray*}
\begin{eqnarray*}
  I_{22}&=&s\lambda^2\int_{-T}^{T} \int_0^1 (a w_x)_x \, \varphi w(| \psi_t|^2-a|\psi_x|^2))\,dx\,dt\\
  &=&-s \lambda^2 \int_{-T}^{T} \int_0^1 \varphi  a | w_x|^2(| \psi_t|^2-a|\psi_x|^2)) \,dx\,dt\\
  &&-\dfrac{s \lambda^2}2 \int_{-T}^{T} \int_0^1 \varphi |w|^2 \left( (|a_x|^2+aa_{xx}) |\psi_x|^2+ 4aa_x\psi_x \psi_{xx} 
  + 2a(a \psi_x)_x \psi_{xx} \right) \,dx\,dt\\
  &&+\dfrac {s \lambda^3}2 \int_{-T}^{T} \int_0^1 \varphi  |w|^2 (a \psi_x)_x (| \psi_t|^2-a|\psi_x|^2))\,dx\,dt\\
  &&+\dfrac {s \lambda^4}2 \int_{-T}^{T} \int_0^1 \varphi  |w|^2 a |\psi_x|^2(| \psi_t|^2- a|\psi_x|^2))\,dx\,dt\\
  &&-s \lambda^3\int_{-T}^{T} \int_{\Omega} \varphi |w|^2a \psi_x\left( a_x |\psi_x|^2 + 2a \psi_x \psi_{xx} \right) \,dx\,dt
\end{eqnarray*}
and
\begin{eqnarray*}
  I_{23} & = & 2s\lambda \int_{-T}^T \int_0^1  (aw_x)_x \, \varphi \left( \psi_t w_t - a\psi_x w_x \right)\, dx\,dt \\
  & = & s \lambda \int_{-T}^{T} \int_0^1 \varphi a |w_x|^2 (\psi_{tt} + a \psi_{xx}) \, dx\,dt \\
  && + s\lambda^2 \int_{-T}^{T} \int_0^1 \varphi a |w_x|^2 \left(  |\psi_t|^2 + a|\psi_x|^2 \right) \, dx\,dt
  - 2s \lambda^2 \int_{-T}^{T}\int_0^1 \varphi a \psi_x \psi_t \, w_x w_t \, dx\,dt \\
  && - s\lambda \int_{-T}^{T}\left[ a(1)^2 |w_x(1,t)|^2 \varphi(1,t)\psi_x(1,t) - a(0)^2 |w_x(0,t)|^2\varphi(0, t)\psi_x(0,t) \right] \,dt.
\end{eqnarray*}

   Finally,
   \[
I_{31} = (\alpha- 1) s^3 \lambda^3 \int_{-T}^{T} \int_0^1 \varphi^3 |w|^2 (| \psi_t|^2 - a |\psi_x|^2)(\psi_{tt} - (a\psi_x)_x) \, dx\,dt,
   \]
   \[
I_{32} = -s^3 \lambda^4 \int_{-T}^{T} \int_0^1 \varphi^3 |w|^2 (| \psi_t|^2 - a |\psi_x|^2)^2 \, dx\,dt
   \]
and
\begin{eqnarray*}
  I_{33} & = & -2s^3\lambda^3 \int_{-T}^T \int_0^1 \varphi^3 w  (| \psi_t|^2 - a |\psi_x|^2) \, (\psi_t w_t - a \psi_x w_x) \, dx\,dt \\
  & = & s^3\lambda^3 \int_{-T}^T \int_0^1 \varphi^3 |w|^2 (|\psi_t|^2 - a |\psi_x|^2) (\psi_{tt} - (a\psi_x)_x) \, dx\,dt \\
  & & + 2s^3 \lambda^3 \int_{-T}^T \int_0^1 \varphi^3 |w|^2 \left( |\psi_t|^2 \psi_{tt}  + aa_x\psi_x |\psi_x|^2 + a^2|\psi_x|^2 \psi_{xx}\right)  \, dx\,dt \\
  && +3 s^3\lambda^4 \int_{-T}^T \int_0^1 \varphi^3 |w|^2  (| \psi_t|^2 - a |\psi_x|^2)^2  \, dx\,dt.
\end{eqnarray*}

Gathering together all terms $I_{ij}$ for $i, \, j \in \{ 1, 2, 3 \}$, we obtain

\begin{align*}
  I &= \int_{-T}^T \int_0^1 (P_1 w) \, (P_2 w)  \, dx\,dt \\
  &= s\lambda \int_{-T}^T \int_0^1 \varphi |w_t|^2 \left( 2\psi_{tt} - \alpha (\psi_{tt} - (a\psi_x)_x) \right)  \, dx\,dt \\
  & + s \lambda \int_{-T}^T \int_0^1 \varphi a |w_x|^2 \left( \alpha (\psi_{tt} - (a\psi_x)_x) + 2(a\psi_x)_x - a_x\psi_x \right) \, dx\,dt \\
  & + 2s\lambda^2 \int_{-T}^T \int_0^1 \varphi \left( |w_t|^2 |\psi_t|^2 - 2 a \psi_x \psi_t w_x w_t + a^2 |w_x|^2 |\psi_x|^2 \right) \, dx\,dt \\
  & + 2s^3\lambda^4 \int_{-T}^T \int_0^1 \varphi^3 |w|^2 \left( |\psi_t|^2 - a |\psi_x|^2 \right)^2  \, dx\,dt \\
  & + s^3\lambda^3 \int_{-T}^T \int_0^1 \varphi^3 |w|^2  (| 2\psi_t|^2 \psi_{tt} + aa_x\psi_x|\psi_x|^2 + 2a^2 |\psi_x|^2\psi_{xx})  \, dx\,dt \\
  & + \alpha s^3 \lambda^3 \int_{-T}^{T} \int_0^1 \varphi^3 |w|^2 (| \psi_t|^2 - a |\psi_x|^2)(\psi_{tt} - (a\psi_x)_x) \, dx\,dt \\
  & - s\lambda \int_{-T}^{T}\left( a(1)^2|w_x(1,t)|^2 \varphi(1,t)\psi_x(1,t) - a(0)^2|w_x(0,t)|^2\varphi(0, t)\psi_x(0,t) \right) \, dx\,dt \\
  & + X_0,
\end{align*}
where $X_0$ is the sum of all ``lower order terms'':
   \[
|X_0| \leq M s\lambda^4 \int_{-T}^T \int_0^1 \varphi^3 |w|^2 \, dx\,dt .
   \]

   Let us analyze the high order terms arising in the previous expression of~$I$. First, remark that
   \begin{equation}\label{eq-geq-0}
s\lambda^2 \int_{-T}^T \int_0^1 \varphi \left( |w_t|^2 |\psi_t|^2 - 2 a \psi_x \psi_t w_x w_t + a^2 |w_x|^2 |\psi_x|^2 \right) \,dx\,dt \geq 0.
   \end{equation}
   
   Secondly, notice that, under the assumption $a \in \mathcal{A}(x_0, a_0)$, if $\beta$ satisfies~\eqref{def_beta}, we can choose $\alpha$ in such a way that the terms of order $s\lambda$ are positive. Indeed, we have in this case
   \begin{equation}
   \nonumber
-a(x) - (x-x_0)a_x(x) < \beta < a(x) + \dfrac{1}{2}(x-x_0)a_x(x) \quad \forall x \in [0,1],
   \end{equation}
whence
   \begin{equation}
   \nonumber
\frac{2\beta}{\beta + a(x) + (x-x_0)a_x(x)} < \frac{2a(x) + (x-x_0)a_x(x)}{\beta + a(x) + (x-x_0)a_x(x)} \quad \forall x \in [0, 1].
   \end{equation}
Let $\alpha$ satisfy
   \begin{equation}
   \nonumber
\sup_{[0,1]} \left( \frac{2\beta}{\beta + a(x) + (x-x_0)a_x(x)} \right) < \alpha < \inf_{[0,1]} \left( \frac{2a(x) + (x-x_0)a_x(x)}{\beta + a(x) + (x-x_0)a_x(x)} \right) .
   \end{equation}
   Then, an explicit computation of the derivatives of $\psi$ shows that
   \[
2\psi_{tt} - \alpha (\psi_{tt} - (a\psi_x)_x) > 0 \ \text{ and } \ \alpha (\psi_{tt} - (a\psi_x)_x) + 2(a\psi_x)_x-a_x\psi_x > 0  \ \text{ in } \ [0,1]\times[-T,T]
   \]
and, consequently,
   \begin{equation}
   \nonumber
\begin{array}{l}
  \dis s\lambda \int_{-T}^T \int_0^1 \varphi |w_t|^2 \left( 2\psi_{tt} - \alpha (\psi_{tt} - (a\psi_x)_x) \right)  \, dx\,dt \\
\noalign{\smallskip}
  \dis \quad \ \ \ \ + \ s\lambda \int_{-T}^T \int_0^1 \varphi a |w_x|^2 \left( \alpha (\psi_{tt} - (a\psi_x)_x) + 2(a\psi_x)_x - a_x\psi_x \right) \, dx\,dt \\
\noalign{\smallskip}
  \dis \ \ \ \ \geq Ms\lambda \int_{-T}^T \int_0^1 \varphi |w_t|^2 \, dx\,dt + Ms\lambda \int_{-T}^T \int_0^1 \varphi |w_x|^2 \, dx\,dt.
\end{array}
   \end{equation}

   The remaining terms in $I$ can be written in the form
   \begin{align*}
  &2s^3\lambda^4 \int_{-T}^T \int_0^1 \varphi^3 |w|^2 \left( |\psi_t|^2 - a |\psi_x|^2 \right)^2  \, dx\,dt \\
  & \ \ \ \ \ + s^3\lambda^3 \int_{-T}^T \int_0^1 \varphi^3 |w|^2  (2|\psi_t|^2 \psi_{tt} + aa_x\psi_x|\psi_x|^2 + 2a^2 |\psi_x|^2\psi_{xx})  \, dx\,dt \\
  & \ \ \ \ \ + \alpha s^3 \lambda^3 \int_{-T}^{T} \int_0^1 \varphi^3 |w|^2 (| \psi_t|^2 - a |\psi_x|^2)(\psi_{tt} - (a\psi_x)_x) \, dx\,dt \\
  & \ \ =s^3 \lambda^3 \int_{-T}^T \int_0^1 \varphi^3 |w|^2 F_\lambda(x,Y(x,t)) \, dx\,dt ,
   \end{align*}
where $Y:= |\psi_t|^2 - a |\psi_x|^2$ and
   \[
\left\{
\begin{array}{l}
\dis F_\lambda(x,Y):= 2\lambda Y^2 + \left( 2\psi_{tt} + \alpha (\psi_{tt} - (a(x)\psi_x)_x)\right) Y \\
\dis \phantom{F_\lambda(x,Y):= } +\ a(x)|\psi_x|^2(2\psi_{tt} + a_x(x) \psi_x +2a(x)\psi_{xx}) \\
\dis \phantom{F_\lambda(x,Y) } = 2\lambda Y^2 + \left(4\beta + \alpha(2\beta + a(x) + (x-x_0)a_x(x)) \right) Y \\
\dis \phantom{F_\lambda(x,Y):= } +\ 8a(x)(x-x_0)^2(-2\beta + 2a(x) +(x-x_0)a_x(x)).
\end{array}
\right.
   \]
Since $F_\lambda$ is polynomial of the second degree in~$Y$, one has
   \[
\begin{array}{l}
\dis F_\lambda(x,Y) \geq 8a(x)(x-x_0)^2(-2\beta + 2a(x) +(x-x_0)a_x(x)) \\
\dis \phantom{F_\lambda(x,Y) \geq } - {1 \over 8\lambda}[4\beta + \alpha(2\beta + a(x) + (x-x_0)a_x(x))]^2
\end{array}
   \]
for all $x \in [0,1]$ and~$Y \in \R$. Therefore, if $\beta$ satisfies~\eqref{def_beta}, for $\lambda$ large enough (depending on $x_0$ and~{\color{black}$\|a\|_{C^3([0,1])}$), } we obtain:
   \begin{equation}\label{eq-geq-2}
s^3 \lambda^3 \int_{-T}^T \int_0^1 \varphi^3 |w|^2 F_\lambda(X) \, dx\,dt \geq M s^3 \lambda^3 \int_{-T}^T \int_0^1 \varphi^3 |w|^2 \, dx\,dt.
   \end{equation}

   Putting together the estimates \eqref{eq-geq-0}--\eqref{eq-geq-2}, the following is found:
   \begin{equation}\label{eq-geq-3}
\begin{array}{l}
\dis \int_{-T}^T \int_0^1 (P_1 w) \, (P_2 w) \, dx\,dt  \geq \ Ms\lambda \int_{-T}^T \int_0^1 \varphi \left( |w_t|^2 + |w_x|^2 \right) \, dx\,dt \\
\noalign{\smallskip}
\dis \qquad +  M s^3 \lambda^3 \int_{-T}^T \int_0^1 \varphi^3 |w|^2 \, dx\,dt \\
\noalign{\smallskip}
\dis \qquad - M s\lambda \int_{-T}^{T} |w_x(1,t)|^2 \, dx\,dt - M s\lambda^4 \int_{-T}^T \int_0^1 \varphi^3 |w|^2 \, dx\,dt.
\end{array}
   \end{equation}   
   
   On the other hand, recalling the definition of $P, \, P_1, \, P_2$ and $R$, we observe that
   \[
\int_{-T}^T \int_0^1 (|P_1 w|^2 \!+\! |P_2 w|^2) \, dx\,dt + 2 \int_{-T}^T \int_0^1 (P_1 w)( P_2 w)\, dx\,dt = \int_{-T}^T \int_0^1 |Pw \!-\! Rw|^2 \, dx\,dt
   \]
It is not difficult to see that there a exists $M$ such that
   \[
\int_{-T}^T \int_0^1 |Pw - Rw|^2 \, dx\,dt \leq M \int_{-T}^T \int_0^1 |Pw |^2 \, dx\,dt + Ms^2\lambda^2 \int_{-T}^T\int_0^1 \varphi^2 |w|^2\, dx\,dt.
   \]
   In particular, we have
   \begin{equation}\label{eq-geq-4}
\begin{array}{l}
\dis \int_{-T}^T \int_0^1 (P_1 w)( P_2 w)\, dx\,dt \leq M \int_{-T}^T \int_0^1 |Pw |^2 \, dx\,dt + Ms^2\lambda^2 \int_{-T}^T\int_0^1 \varphi^2 |w|^2\, dx\,dt
\end{array}
   \end{equation}   
and combining \eqref{eq-geq-3} and~\eqref{eq-geq-3} we obtain:
   \begin{align*}
& s\lambda \int_{-T}^T \int_0^1 \varphi \left( |w_t|^2 + |w_x|^2 \right) \, dx\,dt +  s^3 \lambda^3 \int_{-T}^T  \int_0^1 \varphi^3 |w|^2 \, dx\,dt \\
& \ \ \ \ \leq M \int_{-T}^T \int_0^1 |Pw |^2 \, dx\,dt + M s\lambda \int_{-T}^{T} |w_x(1,t)|^2 \, dt \\
& \quad \ \ \ \ + M s\lambda^4 \int_{-T}^T \int_0^1 \varphi^3 |w|^2 \, dx\,dt+ Ms^2\lambda^2 \int_{-T}^T\int_0^1 \varphi^2 |w|^2\, dx\,dt.
   \end{align*}
 
   Obviously, the last two terms in the right hand side can be absorbed by the second term in the left for $s$ large enough. Therefore, there exists $s_0 > 0$, only depending on $x_0$, $a_0$, {\color{black} $\Vert a\Vert_{C^3([0,1])}$}  and~$T$, such that, for all~$s > s_0$, one has:
   \begin{align}
& s\lambda \int_{-T}^T \int_0^1 \varphi \left( |w_t|^2 + |w_x|^2 \right) \, dx\,dt +  s^3 \lambda^3 \int_{-T}^T  \int_0^1 \varphi^3 |w|^2 \, dx\,dt \nonumber \\
& \ \ \ \ \leq M \int_{-T}^T \int_0^1 |Pw |^2 \, dx\,dt  + M s\lambda \int_{-T}^{T} |w_x(1,t)|^2 \, dt. \label{eq:CarlW}
   \end{align}
Since $w = ve^{s\varphi}$ and $Pw = e^{s\varphi}Lv$, we can easily rewrite \eqref{eq:CarlW} in the form~\eqref{eq:CarlTbis}.

   This ends the proof.
\Fin

\

   In the remaining part of the Appendix, we will use the Carleman estimate \eqref{eq:CarlTbis} to prove Theorem~\ref{th:Carleman}.
   
   Thus, let us assume that \eqref{eq:Tlarge} holds , $w \in L^2(-T, T; H_0^1(0,1))$, $L w \in L^2((0, 1) \times (-T, T))$ and $w_x(1, \cdot) \in L^2(-T, T)$. Thanks to \eqref{eq:Tlarge}, there exists $\eta \in (0, T)$ and $\varepsilon > 1$ such that
   \[
(1 - \varepsilon)(T - \eta)\beta \geq \max_{[0,1]} \, a(x)^{1/2}(x-x_0).
   \]
   
   Moreover, simple computations show that, for every $t \in (-T, -T+\eta) \cup (T-\eta, T)$, the function~$\psi(\cdot\,,t)$ satisfies:
   \begin{equation}\label{eq:condpsi}
\begin{cases}
\dis (1 - \varepsilon) \min_{[0,1]} \,|\psi_t(x,t)| \geq \max_{[0,1]} \,a(x)^{1/2} |\psi_x(x,t)| \\
\noalign{\smallskip}
\dis \max_{[0,1]} \,\psi(x,t) < \min_{[0,1]} \,\psi(x, 0).
\end{cases}
   \end{equation}
Let $\chi \in C_c^\infty(\mathbb{R})$ a cut-off function such that $0 \leq \chi \leq 1$ and
   \[
\chi(t) = \begin{cases}
           1, & \textrm{ if }  |t| \leq T-\eta \\
           0, & \textrm{ if } |t| \geq T
          \end{cases}
   \]
Then we can apply Theorem \ref{th:CarlemanV} to the function $\tilde w:= \chi w$, whence the following Carleman estimate holds
   $$
\begin{array}{c}
\dis s \int_{-T}^T \int_0^1 e^{2s\varphi} \left( |\tilde w_t|^2 + |\tilde w_x|^2 \right) \, dx\,dt +  s^3 \int_{-T}^T  \int_0^1 e^{2s\varphi} |\tilde w|^2 \, dx\,dt \\
\noalign{\smallskip}
\dis \leq M \int_{-T}^T \int_0^1 e^{2s\varphi} |L \tilde w |^2 \, dx\,dt  + M s \int_{-T}^{T} e^{2s\varphi} |\tilde w_x(1,t)|^2 \, dt.
\end{array}
   $$
Since $L \tilde w = \chi Lw + \chi_{tt}w + 2 \chi_t w_t$, we deduce from~\eqref{eq:condpsi} that
   \begin{equation}\label{eq:bef-last}
\begin{array}{l}
\dis s \int_{-T+\eta}^{T-\eta} \int_0^1 e^{2s\varphi} \left( |w_t|^2 + |w_x|^2 \right) \, dx\,dt 
+ s^3 \int_{-T+\eta}^{T-\eta} \int_0^1 e^{2s\varphi} |w|^2 \, dx\,dt \\
\noalign{\smallskip}
\dis \quad \leq M \int_{-T}^T \int_0^1 e^{2s\varphi} |Lw |^2 \, dx\,dt  + M s \int_{-T}^{T} e^{2s\varphi} |w_x(1,t)|^2 \, dt \\
\noalign{\smallskip}
\dis \quad \quad + M\int_{-T}^{-T+\eta} \int_0^1 e^{2 s\varphi} (|w_t|^2 + |w|^2) \, dx\,dt
+ M\int_{T-\eta}^T \int_0^1 e^{2 s\varphi} (|w_t|^2 + |w|^2) \, dx\,dt.
\end{array}
   \end{equation}
   
   Let us denote by $E_s = E_s(t)$ the {\it energy} associated to the operator $L$, that is,
   \begin{equation}\label{eq:energy}
E_s(t):= \frac{1}{2}\int_0^1 e^{2s\varphi} \left( |w_t|^2 + a|w_x|^2 \right) \, dx .
   \end{equation}
   Then, the argument employed in the proof of Theorem~2.5 in~\cite{Bau11} (using the modified energy given by \eqref{eq:energy}) can be used to deduce \eqref{eq:CarlT} from~\eqref{eq:bef-last}.


\end{document}